\documentclass{amsart}

\usepackage{amsfonts, amssymb, graphicx}
\usepackage{amsmath}
\usepackage{amsthm}
\usepackage{color}
\numberwithin{equation}{section}

\usepackage{caption}
\usepackage{subcaption}

\newtheorem{theorem}{Theorem}[section]
\newtheorem{lemma}[theorem]{Lemma}

\newcommand{\e}{\mathrm{e}}

\newcommand{\C}{\mathbb{C}}
\newcommand{\N}{\mathbb{N}}

\numberwithin{equation}{section}

\title[Extension of solutions of first-order difference equations]{On the extension of analytic solutions of first-order difference equations}
\author[R. Halburd]{Rod Halburd}
\address{Rod Halburd \\ Department of Mathematics \\ University College London \\ Gower Street \\ London WC1E 6BT, UK}

\author[R. J. Korhonen]{Risto Korhonen}
\address{Risto Korhonen \\ Department of Physics and Mathematics \\ University of Eastern Finland \\ P.O. Box 111, FI-80101 Joensuu, Finland}
\email{risto.korhonen@uef.fi}

\author[Y. Liu]{Yan Liu}
\address{Yan Liu \\ Department of Physics and Mathematics \\ University of Eastern Finland \\ P.O. Box 111, FI-80101 Joensuu, Finland}
\email{liuyan@student.uef.fi}
\thanks{The second author thanks the support of the China Scholarship Council (No. 202006330038)}

\author[T. Meng]{Techheang Meng}
\address{Techheang Meng \\ Department of Mathematics \\ University College London \\ Gower Street \\ London WC1E 6BT, UK}


\begin{document}

\begin{abstract}
We will consider first-order difference equations of the form 
	\begin{equation*}
	y(z+1) = \frac{\lambda y(z)+a_2(z)y(z)^2
	+\cdots+a_p(z)y(z)^p}{1 + b_1(z)y(z)+\cdots+b_q(z)y(z)^q},
	\end{equation*}
 where $\lambda\in\mathbb{C}\setminus\{0\}$ and the coefficients $a_j(z)$ and $b_k(z)$ are meromorphic.
 When existence of an analytic solution can be proved for large negative values of $\Re(z)$, the equation determines a unique extension to a global meromorphic solution.
 In this paper we prove 
the existence of non-constant meromorphic solutions when the coefficients satisfy $|a_{j}(z)|\leq \nu^{|z|}$ and $|b_{k}(z)|\leq \nu^{|z|}$ for some $\nu<|\lambda|$ in a half-plane.  Furthermore, when a solution exists that is analytic for large positive values of $\Re(z)$, the equation determines a unique extension to a global solution that will generically have algebraic branch points.  We analyse a particular constant coefficient equation, $y(z+1)=\lambda y(z)+y(z)^2$, $0<\lambda<1$, and describe in detail the infinitely-sheeted Riemann surface for such a solution.  We also describe solutions with natural boundaries found by Mahler.
\end{abstract}

\maketitle

\section{Introduction}
In this paper we will study the existence of solutions in the complex domain of first-order difference equations of the form
	\begin{equation}\label{eq1.1}
	y(z+1)=F(z,y(z)),
	\end{equation}
where $F(z,y)$ is a rational function of $y$ with coefficients that are meromorphic functions of $z$.
The existence of meromorphic solutions of the constant coefficient first-order nonlinear difference equation
	\begin{equation}\label{eq1.3}
	y(z+1)=R(y(z)),
	\end{equation}
where $R$ is a rational function with constant coefficients, has a long history going back to Koenigs \cite{koenigs:84}, Poincar\'e \cite{poincare:90}, Fatou \cite{fatou1919,fatou1920a,fatou1920b} and Julia \cite{p2julia}.  Clearly if $y(z)$ is a meromorphic solution of equation \eqref{eq1.1} then so is
\begin{equation}
\label{periodic}
y(z+p(z)),
\end{equation} 
where $p$ is any entire period one function.
It is natural to try to eliminate some of this freedom by restricting the class of solutions of meromorphic solutions of equation \eqref{eq1.1}, for example 
	\begin{equation}\label{eq1.2}
	y(z)\longrightarrow\gamma \quad\textrm{as}\quad \Re(z)\longrightarrow -\infty.
	\end{equation}
This is only possible if $\gamma$ is a fixed point of $R$, i.e. $R(\gamma)=\gamma$.

The main idea to prove global existence is first to show that for any $\sigma>0$, there exists $\rho$ such that there is a one-parameter family of solutions $y$ of equation \eqref{eq1.3} subject to \eqref{eq1.2} that is analytic on
	\begin{equation*}
	{D}(\rho,\sigma)=\left\{z: \Re(z) < -\rho, \ |\Im(z)|<\sigma \right\}.
	\end{equation*}
The equation itself can then be used to continue the solutions as meromorphic functions to the infinite strip $-\sigma<\Im(z)<\sigma$. Finally, since $\sigma$ is arbitrary, it follows that these solutions are meromorphic on $\mathbb{C}$.  Similarly, it is possible to find solutions analytic as $\Re(z)\to +\infty$. In this case, the equation can be used to continue the solution over all $\mathbb{C}$, but generically the function so obtained will have infinitely many branch points and will only be defined on an infinitely sheeted cover of $\mathbb{C}$.

The following lemma is due to G. Julia \cite{p2julia}.
\begin{lemma}
\label{julialemma}
Let $R$ be a rational function of degree greater than one.  Then $R$ has a fixed point $\gamma$ such that $\lambda:=R'(\gamma)$ satisfies either
	\begin{itemize}
	\item[(a)] $|\lambda|>1$ or
	\item[(b)] $\lambda=1$.
	\end{itemize}
\end{lemma}
\noindent
This divides the analysis  into two essentially different cases.

In the case (a), the transformation $y(z)=w(\zeta)$, $\zeta=\lambda^z$, maps equation \eqref{eq1.3} 
to the inverse Schr\"oder equation
\[
	w(\lambda z)=R(w(z)).
\]
For any $\alpha\in\mathbb{C}$, the existence of a unique meromorphic solution satisfying $w(0)=\gamma$ and $w'(0)=\alpha$ goes back to Poincar\'e \cite{poincare:90} and it is implicit in earlier work of Koenigs \cite{koenigs:84}. The existence of a meromorphic solution in a neighbourhood of $0$ can be proved in a straightforward manner and once again the equation itself provides the meromorphic extension ot $\mathbb{C}$.
This case has been studied by S. Shimomura in the case when $R$ is a polynomial \cite{p2Simo} and by N. Yanagihara in the case of rational $R$ \cite{p2Yana}.  We have the following.

\begin{theorem} \label{auto-lambdabig}
Let $R$ be a rational function with a fixed point $\gamma$ such that $\lambda:=R'(\gamma)$ satisfies $|\lambda|>1$. Then for any $\alpha\in\mathbb{C}$, equation \eqref{eq1.3} has a unique meromorphic solution in the complex plane such that for all $\sigma>0$,
	\begin{equation}\label{lb-asympt}
	(y(z)-\gamma)\lambda^{-z} \longrightarrow \alpha, \quad \mbox{for}\  -\sigma<\Im(z)<\sigma\quad\textrm{as}\  \ \Re(z)\longrightarrow-\infty.
	\end{equation}
\end{theorem}

Observe that if there is a fixed point of $R$ such that $0<|R'(\gamma)|<1$, then the same reasoning produces the following result.

\begin{theorem} \label{auto-lambdasmall}
Let $R$ be a rational function with a fixed point $\gamma$ such that $\lambda:=R'(\gamma)$ satisfies $|\lambda|<1$.  Then for any $\alpha\in\mathbb{C}$ and any $\sigma>0$, there is a $\rho>0$ such that equation \eqref{eq1.3} has a unique analytic solution in the set
	\begin{equation*}
	\widehat{D}(\rho,\sigma)=\left\{z: \Re(z) > \rho, \quad |\Im(z)|<\sigma \right\}
	\end{equation*}
satisfying
	\begin{equation}\label{ls-asympt}
	(y(z)-\gamma)\lambda^{-z} \longrightarrow \alpha\quad \textrm{as}\quad \Re(z)\longrightarrow +\infty, \quad z\in\widehat{D}(\rho,\sigma).
	\end{equation}
\end{theorem}

Proving the existence of a meromorphic solution of equation \eqref{eq1.3} corresponding to case (b) of Lemma \ref{julialemma} is more subtle.  Existence was proved by Fatou \cite{fatou1919,fatou1920a,fatou1920b} and by T. Kimura \cite{p2kimura1971,p2kimura1973}.
Baker and Liverpool \cite{p2BakerLiverpool} and Azarina \cite{p2Azarina} later showed that all meromorphic solutions of equation \eqref{eq1.3} can be obtained from the solutions constructed in Shimomura, Yanagihara and Kimura using a transformation of the form \eqref{periodic}.

When the degree of $R$ is one, then either equation \eqref{eq1.3} is linear $y(z+1)=a_0y(z)+a_1$ or it is the discrete Riccati equation $y(z+1)=\{a_0(y(z)+a_1\}/\{y(z)+b_1\}$, which is linearised by the transformation $y(z)=a_0w(z-1)/w(z)$.  So in either case the general solution can be given explicitly in terms of an exponential and a period one function.

The purpose of the present paper is threefold.  First, we will provide a more direct proof than that given by Fatou and Kimura of the existence of a meromorphic solution of the autonomous difference equation \eqref{eq1.3} subject to the asymptotic condition \eqref{eq1.2}, where $R(\gamma)=\gamma$ and $R'(\gamma)=1$. To this end in section \ref{sec3}, we will give a proof of Theorem \ref{th2.1}, building on ideas outlined in Halburd and Korhonen \cite{p2HKphD2006}. Our proof is more direct than those by Fatou and Kimura, who first studied solutions to the Abel equation.  The second purpose is to present a detailed analysis of global solutions of the equation
\[
	y(z+1)=\lambda y(z)+y(z)^2,\qquad 0<\lambda<1.
\]
This equation does not admit any non-constant meromorphic solutions, however, from Theorem \ref{auto-lambdasmall} we see that it has a solution that is meromorphic for large positive $\Re (z)$.  In section \ref{sec7} we will describe the analytic extension of this solution to an infinitely-sheeted cover of the complex plane.
Finally, we will consider meromorphic solutions of the non-autonomous difference equation
	\begin{equation}\label{eq2.1}
	y(z+1) = \frac{\lambda y(z)+a_2(z)y(z)^2+\cdots+a_p(z)y(z)^p}{1 + b_1(z)y(z)+\cdots+b_q(z)y(z)^q}=:F(z,y(z)),
	\end{equation}
where $\lambda\in\C\setminus\{0\}$ and for some $\nu<|\lambda|$,
$a_j(z)$ and $b_k(z)$ are meromorphic functions such that for all $\sigma>0$ there exists a $\rho>0$ such that
$|a_{j}(z)|\leq \nu^{|z|}$ and $|b_{k}(z)|\leq \nu^{|z|}$ for all
\[
	z\in D(\rho,\sigma)=\{z:\Re(z)<-\rho,\ |\Im(z)|<\sigma\}.
\] 
In general, the only fixed point that $F(z,y)$ has a function of $y$ is $y=0$ (i.e., $F(z,0)=0$). Note that $\lambda=F_y(z,0)$. We will have to consider the cases $|\lambda|<1$ and $|\lambda|>1$ separately.

\section{Main results}\label{sec2}

For $\rho>0$, let
\begin{equation*}
    D(\rho)=\left\{z:\Re(z)< -\rho \right\}.
\end{equation*}
The first theorem that we will prove is the following.
\begin{theorem}\label{th2.1}
Let $R$ be a rational function with a fixed point $\gamma\in\mathbb{C}$ such that $R'(\gamma)=1$. Let $m$ be the smallest positive integer such that $R^{(m+1)}(\gamma)\neq 0.$ Given  $\alpha\in\mathbb{C}$, $\delta\in(0,1)$ and $\sigma>0$, there exist constants $\rho>0$ and $\beta\in\mathbb{C}$ and exactly $m$ solutions $y(z)$ of \eqref{eq1.3},  meromorphic in the complex plane, such that
	\begin{equation}\label{eq2.2}
	\frac{2}{R''(\gamma)}\cdot\frac{1}{\gamma-y(z)}=z+\alpha+\beta\log z+W(z)
	\end{equation}
if $m=1,$ and
	\begin{equation}\label{eq2.3}
	-\frac{(m+1)!}{mR^{(m+1)}(\gamma)}\cdot\frac{1}{(y(z)-\gamma)^{m}}=z+\alpha+\beta z^{\frac{m-1}{m}}+W(z)
	\end{equation}
if $m\geq2$ and for all $z$ such that $\Re(z)<-\rho$ and $-\sigma<\Im(z)<\sigma$ we have
\begin{equation*}
	|W(z)|\leq |z|^{-1+\delta}.
	\end{equation*}
\end{theorem}

The next two theorems are generalisations of Theorems \ref{auto-lambdabig} and \ref{auto-lambdasmall} to non-autonomous equations.

\begin{theorem} \label{th2.2}
Let $\lambda\in\mathbb{C}$  and $\nu\in\mathbb{R}$ satisfy $|\lambda|>1$, $0<\nu<|\lambda|$ and let $a_2,\ldots,a_p;b_1,\ldots,b_q$ be meromorphic functions such that,
given $\sigma>0$ there is a $\rho>0$ such that $|a_{j}(z)|\leq \nu^{|z|}$ and $|b_{k}(z)|\leq \nu^{|z|}$ for
\[
	z\in D(\rho,\sigma)=\{z:\Re(z)<-\rho,\ |\Im(z)|<\sigma\}.
\]
Then given $\alpha\in\mathbb{C}$, equation \eqref{eq2.1} has a unique meromorphic solution in the complex plane such that
	\begin{equation}\label{eq2.4}
	y(z)\lambda^{-z} \longrightarrow \alpha, \quad \textrm{as}\quad \Re(z)\longrightarrow-\infty,\quad z\in D(\rho,\sigma).
	\end{equation}
\end{theorem}

\begin{theorem} \label{th2.3}
Let $\lambda\in\mathbb{C}$ and $\nu\in\mathbb{R}$ satisfy $0<\nu<|\lambda|<1$ and let $a_2,\ldots,a_p;b_1,\ldots,b_q$ be meromorphic functions such that,
given $\sigma>0$ there is a $\rho>0$ such that $|a_{j}(z)|\leq \nu^{|z|}$ and $|b_{k}(z)|\leq \nu^{|z|}$ for
\[
	z\in D(\rho,\sigma)=\{z:\Re(z)<-\rho,\ |\Im(z)|<\sigma\}.
\]
Then given $\alpha\in\mathbb{C}$, equation \eqref{eq2.1} has a unique solution analytic in $D(\rho,\sigma)$ such that
	\begin{equation}\label{eq2.5}
	y(z)\lambda^{-z} \longrightarrow \alpha, \quad \textrm{as}\quad \Re(z)\longrightarrow+\infty,\quad z\in D(\rho,\sigma).
	\end{equation}
\end{theorem}

\section{The proof of Theorem \ref{th2.1}}\label{sec3}

We will prove Theorem \ref{th2.1} using Banach's fixed point theorem, including more detail than was provided in the earlier proof outlined by two of us \cite[Theorem 4.2]{p2HKphD2006}. 
The advantage of this approach is that the basic idea is very direct and simple. We do not need to transform the difference equation into the Abel's functional equation, which simplifies the procedure significantly. Banach's theorem also immediately guarantees the uniqueness  of the solution, which is not the case in the reasoning based on normal families.

\begin{theorem}\label{th3.1}
(Banach fixed point theorem). If $T$ is a contraction defined on a complete metric space $X$, then $T$ has a unique fixed point in $X$.
\end{theorem}
We begin by expanding the right hand side of \eqref{eq1.3} as a Taylor series around the fixed point $\gamma$. We obtain
	\begin{equation}\label{eq3.1}
	y(z+1)=\sum_{j=0}^{\infty} \frac{1}{j!}R^{(j)}(\gamma) \left(y(z)-\gamma\right)^j,
	\end{equation}
where $\gamma\in\C$ such that $R(\gamma)\not=\infty$. Substituting $y(z)=w(z)+\gamma$, we then have
	\begin{equation}\label{eq3.2}
	w(z+1)=  w(z) + \sum_{j=2}^{\infty} \frac{1}{j!}R^{(j)}(\gamma) w(z)^j.
	\end{equation}
By applying the transformation $w(z)=1/g(z)$ to equation \eqref{eq3.2}, we obtain
	\begin{equation}\label{eq3.3}
	g(z+1)=\frac{g(z)}{1+ \sum_{j=2}^{\infty} \frac{1}{j!}R^{(j)}(\gamma) g(z)^{-j+1} }=g(z)\left(1-
	\frac{1}{2}R''(\gamma)g(z)^{-1}+\cdots \right).
	\end{equation}
Assume first that $R''(\gamma)\not=0$.  Then
	\begin{equation*}
	\tilde y(z)=-\frac{2g(z)}{R''(\gamma)}
	\end{equation*}
satisfies
	\begin{equation}\label{eq3.4}
	\tilde y(z+1)=F(\tilde y(z)),
	\end{equation}
where
	\begin{equation}\label{eq3.5}
	F(\zeta)= \zeta+1+\sum_{j=1}^{\infty} c_j \zeta^{-j},
	\end{equation}	
for some constants $c_j\in\C$.
So $\widetilde{y}(z)$ is a solution of \eqref{eq3.4} if and only if
	\begin{equation*}
	y(z)=-\frac{2}{R''(\gamma)}\frac{1}{\widetilde{y}(z)}+\gamma
	\end{equation*}
is a solution of the equation \eqref{eq1.3}.

We will next prove that equation \eqref{eq3.4} has an analytic solution in the set $D(\rho)$.  It will be done by showing that there exists a function $W(z)$ such that the following three conditions are satisfied:
	\begin{itemize}
	\item[(1)] $W(z)$ is analytic in the domain $D(\rho)$, where $\rho$ is sufficiently large number.
	\item[(2)]  $|W(z)|\leq |z|^{-1+\delta}$, where $\delta\in (0,1)$ is a fixed constant.
	\item[(3)] $Y(z)=z+\alpha+\beta \log z + W(z)$, where $\alpha,\beta\in\C$, is a solution of \eqref{eq3.4}.
	\end{itemize}

We will use Banach's fixed point theorem to find such a function. To apply the theorem, we need a contraction mapping in a Banach space. For this purpose we define a family $X$ of analytic functions $W$ such that
	\begin{equation}\label{eq3.6}
	|W(z)| \leq |z|^{-1+\delta}
	\end{equation}
holds in $D(\rho)$. Then $X$ is a metric space under the $\sup$-norm
	\begin{equation*}
	\| W\| := \sup_{z\in D(\rho) }\left| W(z)\right|.
	\end{equation*}
Moreover, if a sequence $\{W_n\}$ of $X$ is Cauchy, then $\lim_{n\rightarrow\infty} W_n(z)=W(z)$ exists uniformly in $D(\rho)$. Furthermore, $W$ is analytic in $D(\rho)$ by Weierstrass' theorem on uniformly convergent series of analytic functions. Let finally $\varepsilon >0$. Then there exists $N\in\N$ such that
	\begin{equation}\label{eq3.7}
	|W(z)| \leq |W(z)-W_N(z)| + |W_N(z) | \leq \|W-W_N\| + |W_N(z) | \leq  |z|^{-1+\delta} + \varepsilon.
	\end{equation}
Since \eqref{eq3.7} holds for any $\varepsilon >0$,
	\begin{equation}\label{eq3.8}
	|W(z)| \leq |z|^{-1+\delta}
	\end{equation}
for all $z\in D(\rho)$. Therefore $W\in X$, and so $X$ is a complete metric space.

We will now find a suitable operator in the Banach space $X$. Condition $(3)$ yields
	\begin{equation*}
	\begin{split}
	W(z+1)-W(z) &= -\beta \log \left(1+\frac{1}{z}\right) +\frac{c_1}{z+\alpha+\beta\log z +W(z)}\\
	&\quad +\sum_{j=2}^{\infty}\frac{c_j}{(z+\alpha+\beta\log z +W(z))^j}.\\
	\end{split}
	\end{equation*}
By expanding out the first two terms on the right hand side, we get
	\begin{equation}\label{eq3.9}
	\begin{split}
	&W(z+1)-W(z) =\frac{c_1}{z}-\frac{\beta}{z} +\beta \sum_{j=2}^{\infty} (-1)^j \frac{1}{jz^j} \\	&\quad   +\frac{c_1}{z}\sum_{n=1}^{\infty} \left(-\frac{\alpha+\beta\log z+W(z)}{z}\right)^n+\sum_{j=2}^{\infty}\frac{c_j}{(z+\alpha+\beta\log z +W(z))^j}.\\
	\end{split}
	\end{equation}
We want our operator to satisfy condition \eqref{eq3.8}. For this purpose we need the right hand side of \eqref{eq3.9} to tend to zero sufficiently fast as $|z|$ tends to infinity. Therefore we now
fix $\beta= c_1$ to cancel out the $z^{-1}$ term, and obtain
	\begin{equation*}
	\begin{split}
	W(z+1)-W(z) &= \beta \sum_{j=2}^{\infty} (-1)^j \frac{1}{jz^j}  +\frac{c_1}{z}\sum_{n=1}^{\infty} \left(-\frac{\alpha+\beta\log z+W(z)}{z}\right)^n\\	
 &\quad+\sum_{j=2}^{\infty}\frac{c_j}{(z+\alpha+\beta\log z +W(z))^j}.\\
	\end{split}
	\end{equation*}
Therefore,
	\begin{equation*}
	\begin{split}
	\sum_{k=1}^{m} \left(W(z+1-k)-W(z-k)\right) &=   \sum_{k=1}^{m}\frac{c_1}{z-k}\sum_{n=1}^{\infty} \left(-\frac{\alpha+\beta\log (z-k) +W(z-k)}{z-k}\right)^n\\	 &\quad+\sum_{k=1}^{m}\sum_{j=2}^{\infty}\frac{c_j}{(z-k+\alpha+\beta\log (z-k) +W(z-k))^j}\\
	&\quad + \beta \sum_{k=1}^{m} \sum_{j=2}^{\infty} (-1)^j \frac{1}{j(z-k)^j}.
	\end{split}
	\end{equation*}
Letting $m\longrightarrow\infty$ and using \eqref{eq3.8}, we then have
	\begin{equation*}
	\begin{split}
	W(z) &=  \sum_{k=1}^{\infty}\frac{c_1}{z-k}\sum_{n=1}^{\infty} \left(-\frac{\alpha+\beta\log (z-k) +W(z-k)}{z-k}\right)^n\\	
 &\quad+\sum_{k=1}^{\infty}\sum_{j=2}^{\infty}\frac{c_j}{(z-k+\alpha+\beta\log (z-k) +W(z-k))^j}\\
	&\quad + \beta \sum_{k=1}^{\infty} \sum_{j=2}^{\infty} (-1)^j \frac{1}{j(z-k)^j}.
	\end{split}
	\end{equation*}
This suggests that we define the operator
	\begin{equation}\label{eq3.10}
	\begin{split}
	T[W](z) &=  \sum_{k=1}^{\infty}\frac{c_1}{z-k}\sum_{n=1}^{\infty} \left(-\frac{\alpha+\beta\log (z-k) +W(z-k)}{z-k}\right)^n\\	 
 &\quad+\sum_{k=1}^{\infty}\sum_{j=2}^{\infty}\frac{c_j}{(z-k+\alpha+\beta\log (z-k) +W(z-k))^j}\\
	&\quad + \beta \sum_{k=1}^{\infty} \sum_{j=2}^{\infty} (-1)^j \frac{1}{j(z-k)^j}\\
	&=:T_1\left(W(z)\right)+T_2\left(W(z)\right)+T_3\left(z\right).
	\end{split}
	\end{equation}

We have already shown that $X$ is a complete metric space. We will prove next that $T$ is a mapping from $X$ into $X$. For this it is sufficient to show that
	\begin{equation}\label{eq3.11}
	|T[W](z)| \leq |z|^{-1+\delta}
	\end{equation}
for all $z \in D(\rho)$, where $\delta$ is the constant introduced in the requirement $(2)$. 
We will prove our claim in three parts. First, we look at the function $T_1$. By choosing $\rho$ large enough,
	\begin{equation*}
	|\alpha+\beta\log (z-k) +W(z-k)| \leq |z-k|^{\frac{\delta}{2}}
	\end{equation*}
for all $z\in D(\rho)$. Therefore,
	\begin{equation}\label{eq3.12}
	\begin{split}
	\left|T_1\left(W(z)\right)\right| &=  \left|\sum_{k=1}^{\infty}\frac{c_1}{z-k}\sum_{n=1}^{\infty} \left(-\frac{\alpha+\beta\log (z-k) +W(z-k)}{z-k}\right)^n\right|\\
	&\leq \sum_{k=1}^{\infty}\left|\frac{c_1}{z-k}\right|\sum_{n=1}^{\infty} \left|\frac{1}{z-k}\right|^{n(1-\frac{\delta}{2})}\\
	&=|c_1| \sum_{k=1}^{\infty}\left|\frac{1}{z-k}\right|^{2-\frac{\delta}{2}} \frac{|z-k|^{1-\frac{\delta}{2}}}{|z-k|^{1-\frac{\delta}{2}}-1} \\
	&\leq \frac{2|c_1|}{|z|^{1-\frac{3\delta}{4}}} \sum_{k=1}^{\infty}\left|\frac{1}{z-k}\right|^{1+\frac{\delta}{4}},
	\end{split}
	\end{equation}
when $\rho$ is large enough. Furthermore, 
	\begin{equation*}
	\sum_{k=1}^{\infty}\left|\frac{1}{z-k}\right|^{1+\frac{\delta}{4}}  \leq \sum_{k=1}^{\infty}\frac{1}{k^{1+\frac{\delta}{4}}} =:K,
	\end{equation*}
where the constant $K$ depends only on $\delta$.  So we have a majorising series for $T_1[W]$, so the series converges absolutely and uniformly to an analytic function on $D(\rho)$.
 Finally we also have by \eqref{eq3.12}
	\begin{equation}\label{eq3.13}
	\left|T_1\left(W(z)\right)\right| \leq \frac{2|c_1|K}{|z|^{1-\frac{3\delta}{4}}} \leq \frac{1}{3}|z|^{-1+\delta}
	\end{equation}
for all $z\in D(\rho_1)$, where $\rho_1$ is sufficiently large.

We will now find a majorising series for function $T_2$. Let $\varepsilon>0$. 
We note that
	\begin{equation*}
	|z-k+\alpha+\beta\log (z-k) +W(z-k)| \geq (1-\varepsilon)|z-k|
	\end{equation*}
for all $z\in D(\rho)$ when $\rho$ is large enough. Moreover, since the series \eqref{eq3.5} is convergent,  there is a constant $M>0$ such that $|c_j|\leq M$ for all $j\in\N$. Therefore,
	\begin{equation}\label{eq3.14}
	\begin{split}
	\left|T_2\left(W(z)\right)\right| &=  \left|\sum_{k=1}^{\infty}\sum_{j=2}^{\infty}\frac{c_j}{(z-k+\alpha+\beta\log (z-k) +W(z-k))^j}\right|\\
	&\leq M \sum_{k=1}^{\infty}\sum_{j=2}^{\infty}\left|\frac{1}{(1-\varepsilon)(z-k)}\right|^j\\
	&\leq \frac{M}{(1-\varepsilon)^2} \sum_{k=1}^{\infty} \left|\frac{1}{(z-k)}\right|^2 \frac{(1-\varepsilon)|z-k|}{(1-\varepsilon)|z-k|-1} \\
	&\leq \frac{2M}{(1-\varepsilon)^2} \sum_{k=1}^{\infty} \left|\frac{1}{(z-k)}\right|^2 ,\\
	\end{split}
	\end{equation}
when $\rho$ is large enough. Since
	\begin{equation*}
	 \sum_{k=1}^{\infty} \left|\frac{1}{(z-k)}\right|^2 \leq \frac{1}{|z|^{1-\frac{\delta}{2}}} \sum_{k=1}^{\infty} \frac{1}{k^{1+\frac{\delta}{2}}},
	\end{equation*}
we have a majorising series for $T-[W]$ and by \eqref{eq3.14}
	\begin{equation}\label{eq3.15}
	\left|T_2\left(W(z)\right)\right|  \leq \frac{1}{3}|z|^{-1+\delta}
	\end{equation}
for all $z\in D(\rho_2)$, where $\rho_2$ is sufficiently large.

Similarly,
	\begin{equation}\label{eq3.16}
	\left|T_3\left(z\right)\right|  \leq \frac{1}{3}|z|^{-1+\delta}
	\end{equation}
for all $z\in D(\rho_3)$, where $\rho_3$ is sufficiently large. Combining  \eqref{eq3.10}, \eqref{eq3.13}, \eqref{eq3.15} and \eqref{eq3.16} we conclude that
	\begin{equation*}
	|T[W](z)|\leq |T_1\left(W(z)\right)| + |T_2\left(W(z)\right)| +|T_3\left(z\right)| \leq |z|^{-1+\delta}
	\end{equation*}
for all $z\in D(\rho)$, where $\rho:=\max\{\rho_1,\rho_2,\rho_3\}$. Therefore $T$ maps $X$ into itself.

We will now show that $T$ is a contraction. Assume that $W_1$ and $W_2$ belong to $X$. Then
	\begin{equation}\label{eq3.17}
	\begin{split}
	&T[W_1] (z)-T[W_2](z)\\ 
 &=  \sum_{k=1}^{\infty} \sum_{j=1}^{\infty}  \frac{c_j}{(z-k)^j} \left( \sum_{n=0 }^{\infty} \left(-
	\frac{\alpha+\beta\log (z-k) +W_1(z-k)}{z-k}\right)^n\right)^j\\
	&\quad - \sum_{k=1}^{\infty} \sum_{j=1}^{\infty}  \frac{c_j}{(z-k)^j}\left( \sum_{n=0 }^{\infty} \left(-
	\frac{\alpha+\beta\log (z-k) +W_2(z-k)}{z-k}\right)^n\right)^j\\
	&=\sum_{k=1}^{\infty} \sum_{j=1}^{\infty}  \frac{c_j}{(z-k)^j} \left(1+ \sum_{n=1}^{\infty} \left(-\frac{1}{z-k}\right)^n \sum_{i=0}^n\binom{n}{i} \left(\alpha+\beta\log (z-k)\right)^{n-i}
	W_1(z-k)^i\right)^j\\
	&\quad -\sum_{k=1}^{\infty} \sum_{j=1}^{\infty}  \frac{c_j}{(z-k)^j} \left(1+ \sum_{n=1}^{\infty} \left(-\frac{1}{z-k}\right)^n \sum_{i=0}^n\binom{n}{i} \left(\alpha+\beta\log (z-k)\right)^{n-i}
	W_2(z-k)^i\right)^j.\\
	\end{split}
	\end{equation}
Denoting
	\begin{equation*}
	\mathcal{S}(k,z,W)= \sum_{n=1}^{\infty} \left(-\frac{1}{z-k}\right)^n \sum_{i=0}^n\binom{n}{i} \left(\alpha+\beta\log (z-k)\right)^{n-i}W(z-k)^i,
	\end{equation*}
equation \eqref{eq3.17} takes the form
	\begin{equation}\label{eq3.18}
	\begin{split}
	&T[W_1] (z)-T[W_2](z) = \sum_{k=1}^{\infty} \sum_{j=1}^{\infty}  \frac{c_j}{(z-k)^j} \left(\left(1+\mathcal{S}(k,z,W_1)\right)^j-\left(1+\mathcal{S}(k,z,W_2)^j\right)\right)\\
	&= \sum_{k=1}^{\infty} \sum_{j=1}^{\infty}  \frac{c_j}{(z-k)^j} \sum_{l=1}^{j}\binom{j}{l} \left(\mathcal{S}(k,z,W_1)^l-\mathcal{S}(k,z,W_2)^l\right)\\
	&= \sum_{k=1}^{\infty} \sum_{j=1}^{\infty}  \frac{c_j}{(z-k)^j} \left(\mathcal{S}(k,z,W_1)-\mathcal{S}(k,z,W_2)\right)\sum_{l=1}^{j}\binom{j}{l} \sum_{p=0}^{l-1}\mathcal{S}(k,z,W_1)^{p}\mathcal{S}(k,z,W_2)^{l-1-p}. \\
	\end{split}
	\end{equation}

We will now estimate the modulus of $\mathcal{S}$. For this purpose we note that by Stirling's formula there is a constant $C>0$ not depending on $n$ such that
	\begin{equation}\label{eq3.19}
	\max_{i=1,\ldots,n} \binom{n}{i} \leq C\frac{2^n}{\sqrt{n}}
	\end{equation}
for all $n\in\N$. Also, by choosing sufficiently large $\rho$,
	\begin{equation}\label{eq3.20}
	|\alpha+\beta\log (z-k)| \leq |z-k|^{\delta}
	\end{equation}
for all $z\in D(\rho)$, where $\delta$ is the fixed constant in \eqref{eq3.6}. Since $\frac{n+1}{\sqrt{n}}2^n\leq 4^{n}$ for all $n\in\N$, we have
	\begin{equation}\label{eq3.21}
	\begin{split}
	|\mathcal{S}(k,z,W)| &= \left| \sum_{n=1}^{\infty} \left(-\frac{1}{z-k}\right)^n \sum_{i=0}^n\binom{n}{i} \left(\alpha+\beta\log (z-k)\right)^{n-i}
	W(z-k)^i  \right|\\
	 &\leq C \sum_{n=1}^{\infty} \left|\frac{1}{z-k}\right|^n \sum_{i=0}^n \frac{2^n}{\sqrt{n}} \left|z-k\right|^{\delta(n-i)}\left|z-k\right|^{(-1+\delta)i} \\
	 &\leq C \sum_{n=1}^{\infty} \left|\frac{1}{z-k}\right|^n \sum_{i=0}^n\frac{4^{n}}{n+1}\left|z-k\right|^{\delta n-i}\\
        &=C \sum_{n=1}^{\infty} \left|\frac{1}{z-k}\right|^{n(1-\delta)}\frac{4^{n}}{n+1} \sum_{i=0}^n\left|\frac{1}{z-k}\right|^{i}.
	\end{split}
	\end{equation}
By $\left|z-k\right|>1,$ we have
\begin{equation*}
  \sum_{i=0}^n\left|\frac{1}{z-k}\right|^{i}\leq 1+\frac{n}{\left|z-k\right|}.
\end{equation*}
By choosing $n$ suitably large
\begin{equation*}
   1\leq \frac{n}{\left|z-k\right|}
\end{equation*}
for all $z\in D(\rho),$ and then we can get
\begin{equation}\label{eq3.22}
    \sum_{i=0}^n\left|\frac{1}{z-k}\right|^{i}\leq \frac{2n}{\left|z-k\right|}.
\end{equation}
Substituting \eqref{eq3.22} into \eqref{eq3.21}, we can get
\begin{equation*}
   \begin{split}
         |\mathcal{S}(k,z,W)| &\leq C \sum_{n=1}^{\infty} \left|\frac{1}{z-k}\right|^{n(1-\delta)}\frac{4^{n}}{n+1}\cdot \frac{2n}{\left|z-k\right|}\\
         &\leq 2C \sum_{n=1}^{\infty} 4^{n} \left|\frac{1}{z-k}\right|^{n(1-\delta)+1}\\
         &= 2C \left|\frac{1}{z-k}\right| \sum_{n=1}^{\infty} \left|\frac{4}{(z-k)^{(1-\delta)}}\right|^{n}  \\
	&\leq 5C \left|\frac{1}{z-k}\right|^{2-\delta}
   \end{split}
\end{equation*}
for all $W\in X$. By choosing $\rho$ suitably large
	\begin{equation}\label{eq3.23}
	|\mathcal{S}(k,z,W)| \leq \left|\frac{1}{z-k}\right|
	\end{equation}
for all $z\in D(\rho)$ and for all $W\in X$. Moreover, since
	\begin{equation*}
	\begin{split}
	\left|W_1(z-k)^i-W_2(z-k)^i\right|  &=	\left|\left(W_1(z-k)-W_2(z-k)\right) \sum_{j=0}^{i-1} W_1(z-k)^{i-1-j}W_2(z-k)^j\right| \\
	&\leq  n \left\|W_1-W_2\right\|
	\end{split}
	\end{equation*}
for all $i=1,\ldots,n$, we have, by \eqref{eq3.19}, \eqref{eq3.20} and using the fact that $n\sqrt{n}2^n\leq 4^n$ for all $n\in\N$,
	\begin{equation}\label{eq3.24}
	\begin{split}
	|&\mathcal{S}(k,z,W_1)-\mathcal{S}(k,z,W_2)|\\ 
 &= \left| \sum_{n=1}^{\infty} \left(-\frac{1}{z-k}\right)^n \sum_{i=1}^n\binom{n}{i} \left(\alpha+\beta\log (z-k)\right)^{n-i}\left(W_1(z-k)^i-W_2(z-k)^i\right)\right|\\
	&\leq \sum_{n=1}^{\infty} \left|\frac{1}{z-k}\right|^n \sum_{i=1}^n\binom{n}{i} \left|\alpha+\beta\log (z-k)\right|^{n-i}\left|W_1(z-k)^i-W_2(z-k)^i\right|\\
	&\leq C  \sum_{n=1}^{\infty} \left|\frac{1}{z-k}\right|^{n(1-\delta)+\delta}n^{\frac{3}{2}} 2^{n} \left\|W_1-W_2\right\|.\\
	 &\leq C \left|\frac{1}{z-k}\right|^{\delta} \sum_{n=1}^{\infty} \left|\frac{4}{(z-k)^{(1-\delta)}}\right|^{n}\left\|W_1-W_2\right\|\\
	 &\leq 5C \left|\frac{1}{z-k}\right|   \left\|W_1-W_2\right\|\\
	\end{split}
	\end{equation}
for all $z\in D(\rho)$.

We will now turn our attention back to the operator $T$. By \eqref{eq3.18}, \eqref{eq3.19}, \eqref{eq3.23} and \eqref{eq3.24}
	\begin{equation*}
	\big| T[W_1] (z)-T[W_2](z) \big| \leq5C \sum_{k=1}^{\infty} \sum_{j=1}^{\infty} |c_j| j^{\frac{3}{2}}2^j \left|\frac{1}{z-k}\right|^{j+1}   \left\|W_1-W_2\right\|.
	\end{equation*}
Then, there is a constant $B$ such that
	\begin{equation*}
	5C|c_j| j^{\frac{3}{2}} \leq B 2^j
	\end{equation*}
for all $j\in\N$.  Therefore
	\begin{equation*}
	\begin{split}
	\big| T[W_1] (z)-T[W_2](z) \big| &\leq
	 B\sum_{k=1}^{\infty}\left|\frac{1}{z-k}\right| \sum_{j=1}^{\infty}  \left|\frac{4}{z-k}\right|^{j}   \left\|W_1-W_2\right\|\\
	 & = 4B\sum_{k=1}^{\infty}  \left|\frac{1}{z-k}\right|^2 \frac{|z-k|}{|z-k|-4}   \left\|W_1-W_2\right\|\\
	 &\leq 5B\sum_{k=1}^{\infty}  \left|\frac{1}{z-k}\right|^2  \left\|W_1-W_2\right\|\\
	 \end{split}
	\end{equation*}
for all $z\in D(\rho)$ when $R$ is sufficiently large. Hence
	\begin{equation}\label{eq3.25}
	\begin{split}
	\big| T[W_1] (z)-T[W_2](z) \big|
	&\leq \frac{5B}{|z|^{\frac{1}{2}}} \sum_{k=1}^{\infty} \frac{1}{k^{\frac{3}{2}}}\left\|W_1-W_2\right\| \\
	&= \frac{K}{|z|^{\frac{1}{2}}}\left\|W_1-W_2\right\|,\\
	\end{split}
	\end{equation}
where $K$ is a constant. Let $k<1$. Then, by \eqref{eq3.25},
	\begin{equation*}
	\big| T[W_1] (z)-T[W_2](z) \big| \leq k\left\|W_1-W_2\right\|
	\end{equation*}
for all $z\in D(\rho)$ when $R$ is large enough. Hence
	\begin{equation*}
	\big\| T[W_1] (z)-T[W_2](z) \big\| \leq k\left\|W_1-W_2\right\|,
	\end{equation*}
and so $T:X\longrightarrow X$ is a contraction in $z\in D(\rho)$ as asserted.

By Theorem \ref{th3.1} the mapping $T:X\longrightarrow X$ has a unique fixed point. Hence the existence of a function $W(z)$ satisfying $(1), (2)$ and $(3)$ is proved. In particular,
	\begin{equation*}
	Y(z)=z+\alpha+\beta \log z + W(z)
	\end{equation*}
is a solution of \eqref{eq3.4}. The final step is to continue the analytic solution $Y(z)$ into a meromorphic solution in the whole complex plane by using equation \eqref{eq2.1}.

We will finally outline the proof in the case $R''(\gamma)=0$. Full details are not presented, since the reasoning is very similar to that of the case $R''(\gamma)\not=0$. Assume now that
	\begin{equation*}
	R^{(j)}(\gamma)=0
	\end{equation*}
for all $j=2,\ldots,m$, and that $R^{(m+1)}(\gamma)\not=0$. The transformation
	\begin{equation*}
	\tilde y(z)=-\frac{(m+1)! g(z)^m}{m R^{(m+1)}(\gamma)}
	\end{equation*}
maps equation \eqref{eq3.3} into
	\begin{equation}\label{eq3.26}
	\tilde y(z+1)=\widetilde{F}(\tilde y(z)),
	\end{equation}
where
	\begin{equation}\label{eq3.27}
	\widetilde{F}(z)= z+1+\sum_{j=m}^{\infty} c_j z^{1-\frac{j+1}{m}}.
	\end{equation}	
Let $\delta\in(0,\frac{1}{m})$. The required Banach space $\widetilde{X}$ is the family of analytic functions $W$ such that
	\begin{equation}\label{eq3.28}
	|W(z)| \leq |z|^{-\frac{1}{m}+\delta},
	\end{equation}
holds in $D(\rho)$. By applying Banach's fixed point theorem with a similar reasoning as in the case $R''(\gamma)\not=0$, but instead of \eqref{eq3.10}, using the operator
  \begin{equation*}
   \begin{split}
   \widetilde{T}[W](z) &=  \sum_{k=1}^{\infty}\frac{c_m}{(z-k)^{\frac{1}{m}}}    \sum_{l=1}^{\infty} \kappa_l(m)\left(  \sum_{n=1}^{\infty} \left(-\frac{\alpha+\beta(z-k)^{1-\frac{1}{m}} +W(z-k)}{z-k}\right)^n\right)^l\\	 
 &\quad+\sum_{k=1}^{\infty}\sum_{j=m+1}^{\infty}\frac{c_{j}}{\left(z-k+\alpha+\beta (z-k)^{1-\frac{1}{m}} +W(z-k)\right)^{\frac{j+1}{m}-1}}\\
	&\quad + \beta  \sum_{k=1}^{\infty} (z-k)^{1-\frac{1}{m}} \sum_{j=2}^{\infty} \xi_j(m) \frac{1}{(z-k)^j},\\
	\end{split}
	\end{equation*}
where $\kappa_l(m)$ and $\xi_j(m)$ are constants depending only on $m$. Showing that $\widetilde{T}$ is indeed a contraction mapping from $\widetilde{X}$ into $\widetilde{X}$ is straightforward task by using similar reasoning as before in the case $R''(\gamma)\not=0$ for the operator $T$. Therefore, we have that for any $\alpha\in\C$ there exists a constant $R>0$, and a unique solution $\widetilde{y}(z)$ of \eqref{eq3.26},  meromorphic in the complex plane, such that
	\begin{equation*}
	\widetilde{y}(z)=z+\alpha+\beta z^{\frac{m-1}{m}} +W(z),
	\end{equation*}
where
	\begin{equation*}
	|W(z)|\leq |z|^{-\frac{1}{m}+\delta}
	\end{equation*}
for all $z\in D(\rho)$, and $\beta=\frac{mc_{m}}{m-1}.$

\section{Examples of branched solutions of Theorem \ref{auto-lambdasmall}}
\label{sec7}

Theorem \ref{auto-lambdasmall} guarantees the existence of a unique solution of equation \eqref{eq1.3} which is meromorphic for large $\Re(z)$, satisfying the asymptotic condition
\[
(y(z)-\gamma)\lambda^{-z}\to\alpha,
\qquad
\Re(z)\to+\infty.
\]
In this section, we will consider this solution in the special case
\begin{equation}
\label{norm-ric}
y(z+1)=\lambda y(z)+y(z)^2,\qquad0<\lambda<1,
\end{equation}
subject to the asymptotic condition
\begin{equation}
\label{th2.3eq:1}
y(z)\lambda^{-z}\to\alpha,
\qquad
\Re(z)\to+\infty
\end{equation}
and show that its analytic continuation results in a function defined on a Riemann surface with an infinite number of sheets over $\mathbb{C}$, which we describe in detail.

Before we describe this solution, we remark that equation (\ref{norm-ric}) admits a variety of interesting solutions that can be obtained by analytic continuation of series solutions.
Under the transformation
\[
y(z)=Y(u)-\frac{\lambda}{2},
\qquad
u=\exp\{-2^z\},
\]
equation (\ref{norm-ric})
becomes
\begin{equation}
\label{th2.3eq:2}
Y(u^2)=Y(u)^2+c,
\end{equation}
where
\[
	c=\frac{\lambda}{2}-\frac{\lambda^2}{4}.
\]
This equation was studied by Wedderburn \cite{wedderburn1922mahler}, where he observed that a solution in the case $c=2$ is a generating function for a combinatorial problem arising in the theory of groups.
Equation (\ref{th2.3eq:2}) belongs to a class of functional equations studied by
Poincar\'e and Picard (see references in
Fatou \cite{fatou1919,fatou1920a,fatou1920b}).
Independently, Mahler \cite{mahler1981special} gave a very detailed analysis of solutions of equation (\ref{th2.3eq:2}) that have a simple pole at $u=0$. Notice that this corresponds to a different class of solutions to equation (\ref{norm-ric}) than that described by Theorem \ref{auto-lambdasmall}.
Given $\alpha\in\mathbb{C}$, a simple consequence of Mahler's work \cite{mahler1981special} is the existence of a solution of equation (\ref{norm-ric}) satisfying the asymptotic condition
\[
	\exp\{2^{-z}\}y(z)\to\alpha,\qquad \Re(z)\to+\infty,\quad \beta<\Im(z)<\gamma,
\]
for some $\beta<\gamma$.

In \cite{mahler1981special}, Mahler showed that when $c>0$, equation (\ref{th2.3eq:2}) admits a solution having a simple pole at $u=0$ with the following properties.
\begin{enumerate}
   	 \item When $c\in(0,1/4]$, the solution is analytic on the unit disc $\{u:|u|<1\}$ and continuous on the unit circle $|z|=1$, which is a natural boundary for $Y$ \cite[Th. 7,8]{mahler1981special}.     	
	 \item When $c>1/4$,  the solution can be analytically continued to the unit disc $\{u:|u|<1\}$ with infinitely many quadratic branches where the unit circle is still a natural boundary \cite[Th.3-6]{mahler1981special}. In this case the solution has infinitely many sheets.
\end{enumerate}
Since $\lambda\in(0,2)$, then $c\in(0,1/4]$. Thus, the solution $y(z)$ of equation (\ref{norm-ric}) that corresponds to case 1 is defined on the horizontal strips
\[
	\log_2\left(\frac{4n-1}2\pi\right)<\Im(z)<\log_2\left(\frac{4n+1}2\pi\right),\qquad n\in\mathbb{Z}.
\]
The boundaries of these strips, i.e. the lines $\Im(z)=\log_2\left([2m+1]\pi/2\right)$, $m\in\mathbb{Z}$, are natural boundaries, so the solution cannot be analytically extended.


Equations (\ref{norm-ric}) and (\ref{th2.3eq:2}) can be mapped to the $q$-difference equation
\begin{equation}
\label{th2.3eq:3}
g(qw)=\lambda g(w)+g^2(w),
\end{equation}
where $w=q^z=1/\{\log(1/u)\}$ and
$y(z)=g(w)=Y(u)-(\lambda/2)$. Here if $q=1/2$, then the solution $g(w)$ that corresponds to Mahler's solution is defined on $\Re(w)>1$ where $\Re(w)=1$ is a natural boundary. 
For other choices of $q$, $g(w)$ can be branched and still have natural boundary.

We now return to the solution described in Theorem \ref{auto-lambdasmall}.
In the constant coefficient case, it is slightly easier to study equation (\ref{th2.3eq:3}) than the equivalent form (\ref{norm-ric}).
The solution of equation (\ref{norm-ric}) characterised by 
Theorem \ref{auto-lambdasmall}, i.e., subject to the asymptotic condition (\ref{th2.3eq:1}), corresponds to the solution $g$ of equation (\ref{th2.3eq:3}) of the form
\begin{equation}
\label{th2.3eq:4}
g(w)=g_1w+\sum_{n\geq2}g_nw^n,
\end{equation}
for small $w$.

By substituting (\ref{th2.3eq:4}) into (\ref{th2.3eq:3}), we note that the leading order term yields $q=\lambda$ and $g_1$ is a free parameter,  reflecting the symmetry
$w\mapsto g_1w$ of equation (\ref{th2.3eq:3}).   When $g_1\ne0$, we use the symmetry to set $g_1=-1$, without loss of generality.
Assuming that $g$ admits a solution of the form (\ref{th2.3eq:4}), then we have $q\in(0,1)$ and $g(w)w^{-1}\to g_1$ as $w\to0$, which corresponds to $g(q^z)q^{-z}\to g_1=-1$ as $q^z\to0$. This is equivalent to $y(z)\lambda^{-z}\to g_1=-1$ as $\Re(z)\to+\infty$ confirming that condition (\ref{th2.3eq:1}) holds. In section \ref{sec7}
we will use (\ref{th2.3eq:4}) to show that $y(z)$ is defined on a Riemann surface with infinitely many sheets over $\mathbb{C}$ with only quadratic branches. 
Moreover, with $g_1=-1$, we will show that all branches on each Riemann sheet belong to certain sequences on the positive real axis. We apply arguments similar to those used by Mahler in obtaining the Mahler solution case 2.

We will show that the solutions of 
equation (\ref{norm-ric})
from Theorem \ref{auto-lambdasmall} are indeed branched. By substituting (\ref{th2.3eq:4}) into (\ref{th2.3eq:3}), we obtain the following relations
\begin{equation}
\label{eq:7.2}
b=\lambda,~g_{j}=\frac{1}{\lambda^j-\lambda}\sum_{s=1}^{j-1}g_sg_{j-s},
\end{equation}
where $g_1$ is arbitrary. 

Now we will show that (\ref{th2.3eq:4}) has a non-zero finite radius of convergence by exhibiting a majorising series. 
Observe that for the appropriate branch of the square root function,
\begin{equation}
\label{eq:7.3}
\tilde{g}(x)=\frac{1}{2}-\sqrt{\frac{1}{4}-\tilde{g}_1x}=\sum_{j\geq1}\tilde{g}_nx^n,
\end{equation}
where for all $n>1$, $\tilde{g}_n=\sum_{s=1}^{n-1}\tilde{g}_s\tilde{g}_{n-s}$.
Since this series has a non-zero radius of convergence and $K|g_n|\leq \tilde{g}_n$, where $K=\sup_{j\geq2}(\lambda-\lambda^j)^{-1}$, then $g(w)$ also has a non-zero radius of convergence about $w=0$.

Next we show that $g$ in not entire.  Assume that the radius of convergence of (\ref{th2.3eq:4}) is infinite, then $g(w)$ is an entire solution of (\ref{th2.3eq:3}). If $g$ has a zero at some point $w_1\in\mathbb{C}$, then from equation (\ref{th2.3eq:3}), $g(\lambda^nw_1)=0$ for all $n\in\mathbb{N}$. Since $\lambda^nw_1\to 0$ and $g$ is analytic at $w=0$, it follows that $g(w)\equiv 0$. If $g(w)\not\equiv 0$, then $g$ is a non-constant entire function, so by Picard's Theorem there exists a point $w_2\in\mathbb{C}$ such that $g(w_2)=1-\lambda\ne 0$.  It follows from equation 
(\ref{th2.3eq:3}) that $g(\lambda^nw_2)=1-\lambda\ne 0$ for all $n\in\mathbb{N}$, but $\lambda^nw_2\to 0$, which contradicts the fact that $g$ is analytic and has a zero at $w=0$.

Hence, (\ref{th2.3eq:4}) has a finite radius of convergence, which we denote by $\hat{r}$. This implies that on the circle $|w|=\hat{r}$, there must be a singular point, which we will call $w_0$. Here, $w_0$ cannot be a pole; otherwise (\ref{th2.3eq:3}) shows that $\{\lambda^nw_0\}$ forms a sequence of poles of $g$ converging to $0$, which contradicts the fact that (\ref{th2.3eq:4}) converges for all $w$ in a neighbourhood of $0$. Similarly, $w_0$ cannot be an algebraic pole of $g$. Therefore, $g$ has a regular branch at $w_0$ and (\ref{th2.3eq:3}) shows that $g(w_0)=-\lambda/2$ and $g(\lambda w_0)=-\lambda^2/4$. 

We will show that $w_0=\hat{r}$ is the unique singularity of $g(w)$ defined by equation (\ref{th2.3eq:4}) on the circle $|w|=\hat r$.  Furthermore, we will show that $g$ can be analytically continued to a function with only square root-type branch points.  The resulting Riemann surface $X$ has an infinite number of sheets over $\mathbb{C}$.  All branch points lie over points of the form $\lambda^{1-n}\hat r$, $n\in\mathbb{N}$. We will describe $X$ in terms of a collection of sheets, each of which is of the form $\mathbb{C}\setminus(\lambda^{1-n}\hat r,+\infty)$, for some $n\in\mathbb{N}$, together with the two boundary lines, each lying over $(\lambda^{1-n}\hat r,+\infty)\subset\mathbb{R}$. Each of these boundary lines contains an infinite number of square root branch points and different segments over $((\lambda^{1-m}\hat r,\lambda^{-m}\hat r)$, $m\ge n$, will be common to two sheets. The restriction of the analytic continuation of $g$ to each such sheet will be called a branch.  We will now prove these claims and provide a more detailed description of how the different branches fit together.

Let $h_0$ be the analytic continuation of $g$ (where defined) along all straight lines in $\mathbb{C}$ through the origin.

\vskip 2 mm
\noindent{\underline{Claim 1:} The function $g(w)$ defined by equation (\ref{th2.3eq:4}) can be continued as an analytic function to every point of the circle $|w|=\hat{r}$, except $w=\hat{r}$, where it has a quadratic branch point, i.e. $w_0=\hat{r}$.
    \begin{proof}[Proof of Claim 1]
    Since $g_1=-1$ and $q=\lambda\in(0,1)$, then (\ref{eq:7.2}) shows that $g_n<0$. Also, $g'(0)=g_1=-1<0$ shows that $g$ is locally injective in a domain $D$ containing $0$, and $g(w)$ is decreasing from the value $0$ on $(0,\epsilon)\subset D$ for some small $\epsilon>0$. Assume that $w_0=\hat{r}\e^{\mathrm{i}\nu}\neq\hat{r}$ for some $\nu\in(0,2\pi)$. Then $g(w_0)=-\lambda/2$ and (\ref{th2.3eq:3}) shows that 
    $\{g(\lambda^nw_0):n\in\mathbb{N}\}$ is a sequence of negative numbers converging to $g(0)=0$. 
    We can choose $N$ sufficiently large such that $g(\lambda^Nw_0)=g(v_0)$ for some point $v_0\in(0,\epsilon)\subset D$ (the point $v_0$ exists by the intermediate value theorem). 
    Using the injectivity of $g$ in $D$, we find that $\lambda^Nw_0=v_0$, which shows that $w_0$ is real and greater than zero, so $w_0=\hat{r}$. 
    \end{proof}
\noindent{\underline{Claim 2:}  The function $g(w)$ can be continued as an analytic function to $\hat X_0:=\mathbb{C}\setminus[\hat r,+\infty)$. Further continuation to the two boundary lines $[\hat r,\infty)^+$ and $[\hat r,\infty)^+$ results in a function $h_0(w)$ analytic at all points of $X_0:=\hat X_0\cup[\hat r,+\infty)^+\cup[\hat r,+\infty)^-$, except at points lying above
 $B:=\{\lambda^{1-n}\hat{r}:n\in\mathbb{N}\}$, each of which is a square root branch point.
    \begin{proof}[Proof of Claim 2]
    Since $w_0=\hat{r}$ is a branch point of $h_0$, then (\ref{th2.3eq:3}) shows that all points lying above $B$ are also branch points.
    Assume that there exists a point $\tilde{w}_0\notin B$ where $h_0$ is branched. Then, using (\ref{th2.3eq:3}), we find that $\lambda^n\tilde{w}_0$ for all sufficiently large $n$ are regular points of $g$. Thus, there exists $N\in\mathbb{N}$ such that $g(\lambda^N\tilde{w}_0)=-\lambda/2,g(\lambda^{N+1}\tilde{w}_0)=-\lambda^2/4$. 

Let $\tilde w_1=\lambda^N\tilde w_0$.  Then $g(\lambda\tilde{w}_1)=-\lambda^2/4$.
Therefore, $g(\lambda^n\tilde{w}_1)=g(\lambda^n\hat{r})$ for all $n\in\mathbb{N}$.  For all $n>0$, $\lambda^n\hat{r}\in D$, and by choosing $n$ sufficiently large, $\lambda^n\tilde{w}_1\in D$.  The injectivity of $g$ on $D$ implies that $\lambda^n\tilde{w}_1=\lambda^n\hat{r}$ for large enough $n$ i.e. $\tilde{w}_1=\hat{r}=w_0$, and so $\tilde w_0\in B$.
    \end{proof}

Now we show that all the branch points of $g$ are quadratic. Let $d_0=\lambda^2/4-\lambda/2$ and since $0<\lambda<1$, then $d_0<0$. We define $d_n=d\pm\sqrt{d_{n-1}}$; then $d_N=0$ for some $N\in\mathbb{N}$ if and only if $\lambda=0$. To see this, we note that $\Im(\sqrt{d_0})\neq0$, so $\Im(d_1)\neq0$. If $\Im(d_{n-1})\neq0$, then $\Im(d_n)=\Im(\sqrt{d_{n-1}})\neq0$ and thus $d_n\neq0,n\in\mathbb{N}$. About the point $w_{-1}:=\lambda\hat{r}$, we have the expansion
\begin{equation}
\label{eq:7.5}
g(w)=-\lambda^2/4+\phi_{0,1}(w-w_{-1})+\sum_{n\geq2}\phi_{0,n}(w-\lambda\hat{r})^n,
\end{equation}
where $\phi_{0,1}=g'(w_{-1})<0$. Moreover, $\phi_{0,n}=n!g^{(n)}(w_{-1})<0$ follows from $g_n<0$. Now, about $w_0=\hat{r}$, using (\ref{eq:7.5}) and (\ref{th2.3eq:3}), we obtain
\begin{equation}
\label{eq:7.6}
\begin{split}
g(w)&=-\lambda/2+\sqrt{\phi_{0,1}\lambda(w-w_{0})+O((w-w_{0})^2)}
\\
&=-\lambda/2+\phi_{1,1}(w-w_0)^{1/2}+\sum_{n\geq2}\phi_{1,n}(w-w_{0})^{n/2},
\end{split}
\end{equation}
where $\phi_{1,1}=\sqrt{g'(w_{-1})\lambda}=\pm \mathrm{i}\sqrt{-g'(w_{-1})\lambda}\neq0$ with two branch choices and $\phi_{1,n}$ can be expressed in terms of $\phi_{1,1}$. Now we fix a branch of $\phi_{1,1}$. Using (\ref{eq:7.6}) and (\ref{th2.3eq:3}) again, about $w_1:=\hat{r}\lambda^{-1}$ we obtain 
\begin{equation}
\label{eq:7.7}
g(w)=-\lambda/2+\sqrt{d_0}\phi_{2,1}(w-w_1)^{1/2}+\sum_{n\geq2}\phi_{2,n}(w-w_{1})^{n/2},
\end{equation}
where $\sqrt{d_0}\phi_{2,1}\neq0$ with two branch choices and $\phi_{2,n}$ can be expressed in terms of $\phi_{2,1}$. Repeating the arguments, we obtain about $w_m:=\hat{r}b^{-m}$,
\begin{equation}
\label{eq:7.8}
g(w)=-\lambda/2+\sqrt{d_{m-2}}\phi_{m,1}(w-w_m)^{1/2}+\sum_{n\geq2}\phi_{2,n}(w-w_{m})^{n/2},
\end{equation}
where $d_{m-2}\phi_{m,1}\neq0$. Hence all branch points of $g$ are quadratic. 

We now specify how we handle the branch cuts at the branch points. On the branch $h_0$ and at the branch point $\hat{r}\lambda^{1-n}$, we take the branch cut along the interval $[\hat{r}\lambda^{1-n},\hat{r}\lambda^{-n})$. Below, we will prove that on each interval $[\hat{r}\lambda^{1-n},\hat{r}\lambda^{-n})$, $h_0$ will connect to different sheets obtained from each branch point $\hat{r}\lambda^{1-n}$ as $n$ varies. We will also show that this property holds for the other Riemann sheets of $g$. Once we fully describe how each sheet connects with the others, we will have a complete picture of the Riemann surface $X$.

We have so far shown that the analytic continuation of $g$ has only quadratic branch points, and $B$ is the set of branch points of $h_0$. However, we do not immediately know the distribution of the branch points on the other sheets. We now investigate the branch points on each sheet, starting with the two sheets we obtain from the branch point $w_0=\hat{r}$. First of all, (\ref{eq:7.7}) has the following form
\begin{equation}
\label{eq:7.9}
\begin{aligned}
g(w)=&-\lambda/2\pm\sqrt{g(\lambda w)+\lambda^2/4}\\
=&-\lambda/2+\sqrt{\phi_{0,1}\lambda(w-w_{0})}(1+\sum_{n\geq2}\lambda^{n-1}\phi_{0,n}\phi_{0,1}^{-1}(w-w_0)^n)^{1/2}\\
=&-\lambda/2\pm\mathrm{i}\sqrt{-g'(w_{-1})\lambda(w-w_{0})}\sum_{n\geq0}\psi_{1,n}(w-w_0)^n,
\end{aligned}
\end{equation}
where $\psi_{1,n}\in\mathbb{R}$ and we fix $\Im(\sqrt{w-w_0})>0$ for $w\in[0,\hat{r}]$. We denote the other branch connected to $h_0$ by the branch point $\hat r$ by $h_1$. Since $h_0$ is the analytic continuation of (\ref{th2.3eq:4}), we must have $\phi_{1,1}=-\mathrm{i}\sqrt{-g'(w_{-1})\lambda}$. To see this, we use (\ref{eq:7.9}) and the fact that $g'(w)<0$ for $w\in[0,\hat{r}]$. From (\ref{eq:7.9}), we obtain
\begin{equation}
\label{eq:7.10}
h_0+h_1=-\lambda.
\end{equation}
This shows that the set of branch points of $h_1$ is also $B$. Here the sheets $h_0$ and $h_1$ connect at the quadratic branch point $\hat{r}$. 
Let $\lambda w=\lambda\hat{r}+t\e^{\mathrm{i}\theta},\theta\in[0,2\pi]$ and $0<t<\lambda(\hat{r}\lambda^{-1}-\hat{r})$. Then on $\lambda\hat{r}+t\e^{\mathrm{i}\theta}$, we obtain 
\begin{equation}
\label{eq:7.12}
h_0(\lambda w)=-\lambda^2/4+\gamma t\e^{\mathrm{i}\theta}+O(t^2),\gamma=h_0'(\lambda\hat{r})<0.
\end{equation}
On $w=\hat{r}+t\e^{\mathrm{i}\theta}\lambda^{-1}$, we obtain
\begin{equation}
\label{eq:7.13}
g(w)=-\lambda/2\pm(\gamma t/\lambda)^{1/2}\e^{\mathrm{i}\theta}+O(t^{3/2}).    
\end{equation}
From the above assumption, we have 
\begin{equation}
\label{eq:7.14}
\begin{aligned}
h_0(w)=&-\lambda/2-\mathrm{i}(-\gamma t)^{1/2}\e^{\mathrm{i}\theta/2}+O(t^{3/2}),\\
h_1(w)=&-\lambda/2+\mathrm{i}(-\gamma t)^{1/2}\e^{\mathrm{i}\theta/2}+O(t^{3/2}).
\end{aligned}
\end{equation}
Here, we observe that if we start on $h_0(w)$, then we will end up with $h_1(w)$ as $\theta$ changes from $0$ to $2\pi$, and return back to the value $h_0(w)$ as $\theta$ changes from $2\pi$ to $4\pi$. We include a figure (Figure \ref{fig:1}) regarding how the sheets $h_0,h_1$ connect at $\hat{r}$. In the second step, we will see that when $t>\lambda(\hat{r}\lambda^{-1}-\hat{r})$, analytically continuing on the circle $w=\hat{r}+t\e^{\mathrm{i}\theta}\lambda^{-1}$ twice does not return back to the initial value (or sheet) where we started. This is due to the fact that on $[\hat{r}b^{-n},\hat{r}b^{-n-1})$, the sheets $h_0,h_1$ connect to new, distinct sheets (see Figure \ref{fig:2}). We would like to remark that $h_0,h_1$ satisfy
\begin{equation}
\label{eq:7.11}
h_0(\lambda w)=h_0^2(w)+\lambda h_0(w)=h_1^2(w)+\lambda h_1(w).
\end{equation}
The functional equation for $h_1$ is actually 
\[
h_1(\lambda w)=-h_1^2(w)-\lambda h_1(w)-\lambda,
\]
which is different from (\ref{th2.3eq:3}). 

Now we proceed to the second step, where there are four sheets resulting from the branch point $\lambda^{-1}\hat{r}$. 
Let $\lambda w=\lambda\hat{r}+t\e^{\mathrm{i}\theta},\theta\in[0,2\pi]$ and $0<t<\lambda(\hat{r}\lambda^{-1}-\hat{r})$. Then on $\lambda\hat{r}+t\e^{\mathrm{i}\theta}$, we obtain (\ref{eq:7.12}) and (\ref{eq:7.14}).
On $w=\hat{r}/\lambda+t\e^{\mathrm{i}\theta}/\lambda^2$, we obtain the expansions of the four branches of $g(w)$ about $\hat{r}/\lambda$, which are 
\begin{equation}
\label{eq:7.15}
\begin{aligned}
&-\lambda/2\pm\sqrt{h_0(\lambda w)+\lambda^2/4} & &=-\lambda/2\pm\sqrt{\lambda^2/4-\lambda/2-(\gamma t)^{1/2}\e^{\mathrm{i}\theta/2}+O(t^{3/2})}\\
& & &=-\frac{\lambda}{2}\pm\left\{\sqrt{\frac{\lambda^2}{4}-\frac{\lambda}{2}}+\left(\frac{\gamma t}{\lambda^2/4-\lambda/2}\right)^{\frac{1}{2}}\e^{\mathrm{i}\theta/2}+O(t^{3/2})\right\},\\
&-\lambda/2\pm\sqrt{h_1(\lambda w)+\lambda^2/4} & &=-\lambda/2\pm\sqrt{\lambda^2/4-\lambda/2+(\gamma t)^{1/2}\e^{\mathrm{i}\theta/2}+O(t^{3/2})}\\
& & &=-\frac{\lambda}{2}\pm\left\{\sqrt{\frac{\lambda^2}{4}-\frac{\lambda}{2}}-\left(\frac{\gamma t}{\lambda^2/4-\lambda/2}\right)^{\frac{1}{2}}\e^{\mathrm{i}\theta/2}+O(t^{3/2})\right\}.
\end{aligned}
\end{equation}
By choosing suitable signs for the branches, we obtain
\begin{equation}
\label{eq:7.16}
\begin{aligned}
&h_0(w)=-\lambda/2-\left\{\sqrt{\lambda^2/4-\lambda/2}+\left(\frac{\gamma t}{\lambda^2/4-\lambda/2}\right)^{1/2}\e^{\mathrm{i}\theta/2}+O(t^{3/2})\right\},\\    
&h_1(w)=-\lambda/2+\left\{\sqrt{\lambda^2/4-\lambda/2}+\left(\frac{\gamma t}{\lambda^2/4-\lambda/2}\right)^{1/2}\e^{\mathrm{i}\theta/2}+O(t^{3/2})\right\}.
\end{aligned}
\end{equation}
Here, we do not restrict $\sqrt{\lambda^2/4-\lambda/2}$ to be a principal branch. Hence, the identities in (\ref{eq:7.16}) are merely symbolic and we may need numerical result to demonstrate which branch of the square root we should choose.

We denote the other two branches corresponding to $-\lambda/2\pm\sqrt{h_1(\lambda w)+\lambda^2/4}$ by 
\begin{equation}
\label{eq:7.17}
\begin{aligned}
&h_{10}(w)=-\lambda/2-\left\{\sqrt{\lambda^2/4-\lambda/2}-\left(\frac{\gamma t}{\lambda^2/4-\lambda/2}\right)^{1/2}\e^{\mathrm{i}\theta/2}+O(t^{3/2})\right\},\\    
&h_{11}(w)=-\lambda/2+\left\{\sqrt{\lambda^2/4-\lambda/2}-\left(\frac{\gamma t}{\lambda^2/4-\lambda/2}\right)^{1/2}\e^{\mathrm{i}\theta/2}+O(t^{3/2})\right\}.    
\end{aligned}
\end{equation}
Here we also obtain
\begin{equation}
\label{eq:7.18}
h_{10}+h_{11}=-\lambda
\end{equation}
and the relations between $h_{10},h_{11}$ and $h_1$ are given by
\begin{equation}
\label{eq:7.19}
h_{1}(\lambda w)=h_{10}^2(w)+\lambda h_{10}(w)=h_{11}^2(w)+\lambda h_{11}(w).
\end{equation}
By observing the value of $\theta\in[0,4\pi]$ in (\ref{eq:7.16}) and (\ref{eq:7.17}), we see that $h_{0}$ and $h_{10}$ connect at $\lambda^{-1}\hat{r}$. For instance, if we start at $h_0(w)$, we will end up with $h_{10}(w)$ as $\theta$ changes from $0$ to $2\pi$, and return back to $h_0(w)$ as $\theta$ changes from $2\pi$ to $4\pi$.
Similarly $h_{1},$ and $h_{11}$ connect at $\lambda^{-1}\hat{r}$ and analytically continuing on the circle $w=\hat{r}/\lambda+t\e^{\mathrm{i}\theta}/\lambda^2$ twice will return to the sheet we started on. We include a diagram of the above details in Figure \ref{fig:2}, showing how these four sheets $h_0,h_1,h_{10}$ and $h_{11}$ connect at $\lambda^{-1}\hat{r}$.
The functional equations of $h_{10}$ and $h_{11}$ can be obtained as follows. Substituting (\ref{eq:7.19}) into $h_1(\lambda w)=-h_1^2(w)-\lambda h_1(w)-\lambda$ yields
\begin{equation}
\label{eq:7.20}
(h_{10}(w)+\lambda/2)^2+(h_1(w)+\lambda/2)^2=\lambda^2/2-\lambda.
\end{equation}
Changing $w\to \lambda w$ in (\ref{eq:7.20}) and eliminating $h_1$ using (\ref{eq:7.20}) yields
\begin{equation}
\label{eq:7.21}
(h_{10}(\lambda w)+\lambda/2)^2+(h_{10}^2(w)+\lambda h_{10}(w)+\lambda/2)^2=\lambda^2/2-\lambda,
\end{equation}
which is the functional equation of $h_{10}$. Using (\ref{eq:7.18}), we can also obtain the functional equation of $h_{11}$. Both equations are different from (\ref{th2.3eq:3}).

Now we show that, on the sheets $h_{10},h_{11}$, the set of branch points is $\{\lambda^{-n}w_0,n\in\mathbb{N}\}$. 
Assume that there is a point $\tilde{w}_0\notin\{\lambda^{-n}\hat{r}\}\cup\{0\}$ which is a branch of $h_{10}$ (here the branches of $h_{11},h_{10}$ are the same by using (\ref{eq:7.18})). Then (\ref{eq:7.19}) shows that either $h_1$ is branched at $\lambda\tilde{w}_0$ or $h_1$ is regular there and $h_1(\lambda\tilde{w}_0)=-\lambda^2/4$. First we show that
\vskip 2mm
\noindent{\underline{Claim 3:}  $h_0(v_0)=h_1(v_1)$ if and only if $v_0=v_1=\hat{r},h_0(\hat{r})=h_1(\hat{r})=-\lambda/2$ i.e. only at the branch point $\hat{r}$.
    \begin{proof}[Proof of Claim 3]
    Assume the contrary that $h_0(v_0)=h_1(v_1)$ for some distinct $v_0$ and $v_1$. Then, by using (\ref{eq:7.11}), we have $h_0(\lambda v_0)=h_0(\lambda v_1)$. Again we obtain $h_0(\lambda^nv_0)=h_0(\lambda^nv_1)$ which yields by injectivity of $D$ that $v_0=v_1$, a contradiction. Therefore, $v_0=v_1$ and $h_0(v_0)=h_1(v_0)=-\lambda/2$. We have already shown that $h_0(\hat{r})=-\lambda/2$. Thus, $v_0=v_1=\hat{r}$.    
    \end{proof}
If $h_1$ is branched at $\lambda\tilde{w}_0$, then we obtain $\lambda\tilde{w}_0\in B=\{\lambda^{1-n}w_0\}$ and so $\tilde{w}_0\in \{\lambda^{-n}w_0\}$ which is a contradiction. So $h_1$ is regular at $\lambda\tilde{w}_0$ and $h_1(\lambda\tilde{w}_0)=-\lambda^2/4$. Since $h_0(\lambda w_0)=-\lambda^2/4$, then $h_1(\lambda\tilde{w}_0)=h_0(\lambda w_0)=-\lambda^2/4$ which contradicts Claim 3. 
Now we show that $0$ cannot be a branch point of $h_{10},h_{11}$. Actually, we can prove a general statement: 
\vskip 2 mm
\noindent{\underline{Claim 4:}  The point $0$ cannot be a branch point on any sheet.
    \begin{proof}[Proof of Claim 4:]
    A new Riemann sheet $H_1$ is obtained from a specific known Riemann sheet $H_0$ through the relation (see (\ref{eq:7.11}) and (\ref{eq:7.18}))
    \begin{equation}
    \label{0branch}    
    H_0(\lambda w)+\lambda^2/4=(H_1(w)+\lambda/2)^2,
    \end{equation}
    where there is a branch point at $\tilde{w}_0$ when $H_0(\lambda \tilde{w}_0)=-\lambda^2/4$. If there exists a Riemann sheet $H_n$, where $0$ is a branch point, we can use (\ref{0branch}) to trace back to the sheet $H_i$ where $H_i(0)=-\lambda^2/4$, which is where the appearance of a branch point would begin. However, we will show that this is not possible. Let $e_0=\lambda^2/4-\lambda<0,e_1=-\lambda/2\pm\sqrt{e_0}$ and $e_n=-\lambda/2\pm\sqrt{\lambda^2/4+e_{n-1}}$. Then, $h_{10}(0),h_{11}(0)=e_1$. Inductively, we can show that for any new Riemann sheet $H$, we have $H(0)=e_j$ for some $j$. Using the same reasoning as when we showed that $d_N\neq0$, we can conclude that $e_j+\lambda^2/4\neq0$. Therefore, $0$ cannot be a branch point of any sheet.    
    \end{proof}
Hence all branch points of $h_{10},h_{11}$ are $\lambda^{-n}w_0,n\in\mathbb{N}$.

Now we describe the situation where we have 8 Riemann sheets in the third step. The 8 sheets are obtained at $\lambda^{-2}\hat{r}$ by 
\begin{equation}
\label{eq:7.22}
\begin{aligned}
&-\lambda/2\pm\sqrt{h_0(\lambda w)+\lambda^2/4},-\lambda/2\pm\sqrt{h_1(\lambda w)+\lambda^2/4},\\
&-\lambda/2\pm\sqrt{h_{10}(\lambda w)+\lambda^2/4},-\lambda/2\pm\sqrt{h_{11}(\lambda w)+\lambda^2/4}.    
\end{aligned}
\end{equation}
The four new sheets come from $-\lambda/2\pm\sqrt{h_{10}(\lambda w)+\lambda^2/4},-\lambda/2\pm\sqrt{h_{11}(\lambda w)+\lambda^2/4}$. We denote them by $h_{100},h_{101},h_{110},h_{111}$ as follows
\begin{subequations}
\label{eq:7.23}
\begin{align}
&h_{10}(\lambda w)=h_{100}(w)^2+\lambda h_{100}(w)=h_{101}(w)^2+\lambda h_{101}(w),\\  &h_{11}(\lambda w)=h_{110}(w)^2+\lambda h_{110}(w)=h_{111}(w)^2+\lambda h_{111}(w).  
\end{align}
\end{subequations}
By a similar argument, we also have
\begin{equation}
\label{eq:7.24}
h_{100}+h_{101}=h_{110}+h_{111}=-\lambda.
\end{equation}
Now we prove that the set of branch points of $h_{100},h_{101},h_{110},h_{111}$ is $\{\lambda^{-n-1}w_0,n\in\mathbb{N}\}$ (we use Claim 4 to conclude that $0$ cannot be a branch point). The proof is similar to the second step, where we assume that $h_{100}$ has a branch at $\tilde{w}_0\notin\{\lambda^{-n-1}w_0\}$ (similar proof for $h_{101},h_{110},h_{111}$). Similar arguments show that $h_{10}$ cannot have a branch point at $\lambda\tilde{w}_0$; otherwise $\tilde{w}_0\in\{\lambda^{-n}w_0\}$ leading to $h_{10}(\lambda\tilde{w}_0)=-\lambda^2/4$. Thus, we again use (\ref{eq:7.19}) ,(\ref{eq:7.11}) and the fact that $h_0(\lambda w_0)=-\lambda^2/4$ to obtain $h_1(\lambda^2\tilde{w}_0)=h_0(\lambda^2w_0)$ which is a contradiction. Hence, all the branch points of $h_{100}$ are $\lambda^{-n-1}w_0,n\in\mathbb{N}$. To see how $h_0,h_1,h_{10},h_{11},h_{100},h_{101},h_{110},h_{111}$ connect at $\lambda^{-2}w_0$, we perform the same analysis to obtain the analogues of (\ref{eq:7.15}) at $\lambda^{-2}w_0$. We include a diagram (see Figure \ref{fig:3}) regarding how the sheets $h_0,h_1,\dots,h_{111}$ connect at $\hat{r}/\lambda^2$.

\begin{figure}[h]
     \centering
     \hspace{-3.4cm}
     \begin{subfigure}[b]{0.3\textwidth}
         \centering
         \includegraphics[scale=1.0]{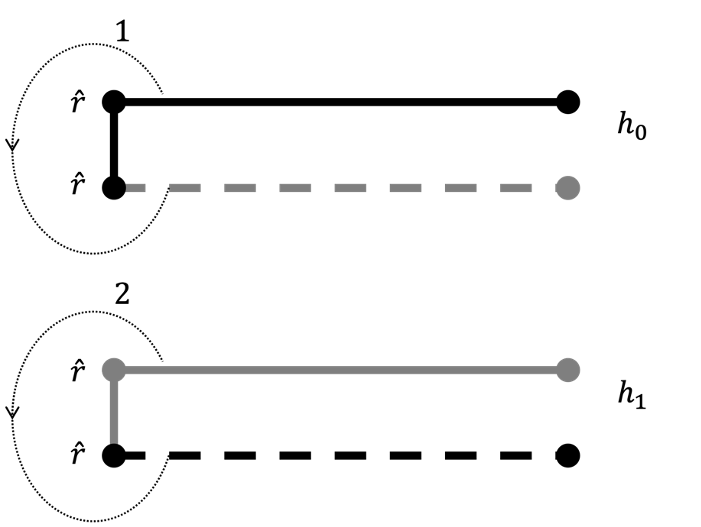}
         \caption{Branch point at $\hat{r}$}
         \label{fig:1}
     \end{subfigure}
     \hspace{3cm}
     \begin{subfigure}[b]{0.3\textwidth}
         \centering
         \includegraphics[scale=1.4]{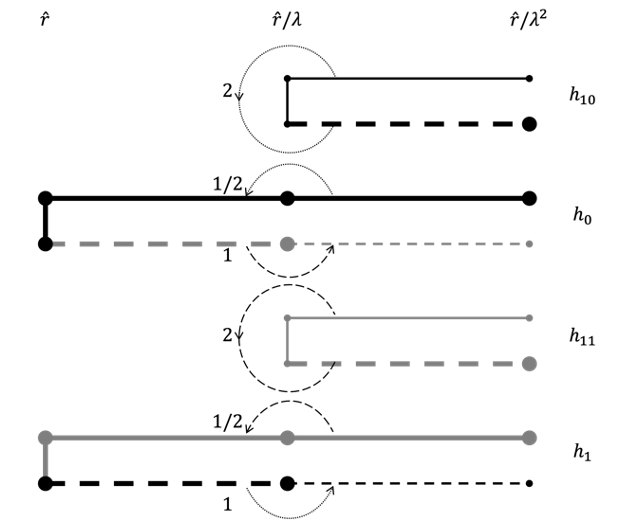}
         \caption{Branch point at $\hat{r}/\lambda$}
         \label{fig:2}
     \end{subfigure}
     \caption{Branch points at $\hat{r}$ and $\hat{r}/\lambda$}
     We cut each sheet along the positive real axis. The round curves represent circles centred at $\hat{r},\hat{r}/\lambda$. Different patterns on the circles represent different circles. The numbers $1/2,1,$ and $2$ denote the number of rounds we trace along the circle.
     \label{fig:0}
\end{figure}

\begin{figure}[h]
    \centering
    \includegraphics[scale=2.0]{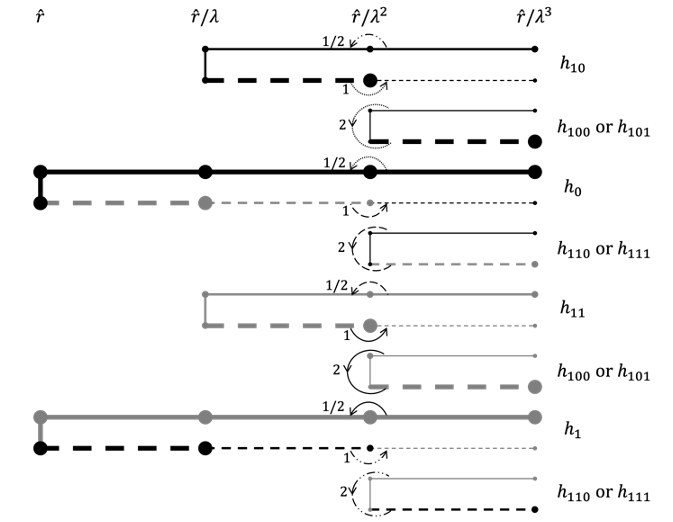}
    \caption{Branch point at $\hat{r}/\lambda^2$}
    \label{fig:3}
    (See Figure \ref{fig:0} for description of the picture)
\end{figure}

By repeating the arguments at the $m$-th step, we obtain, at the quadratic branch $\lambda^{1-m}w_0$, the $2^{m}$ sheets where $2^{m-1}$ are the new sheets. All the branch points of the new sheets are $\lambda^{2-m-n}w_0,n\in\mathbb{N}$. 
When $g_1\neq0,g_1\notin\mathbb{R}_-$, then above results still hold by simply changing $w_0=\hat{r}$ to $w_0=\hat{r}\e^{\mathrm{i}\nu}$ for some $\nu\in[0,2\pi)$. 

Now we return to the domain of $y(z)$. Since $w=b^z=\lambda^z$, the circle $C(0,\hat{r})=\{w:|w|=\hat{r}\}$ is mapped to the vertical line $\Re(z)=\ln\hat{r}/\ln \lambda$\footnote{The value $\hat{r}$ depends on $g_1$ and $\ln\hat{r}/\ln \lambda$ can be negative, compare to the condition in $\hat{D}(\rho,\sigma)$ in section \ref{sec2}.}. The pre-images of $\hat{r}$ are $z=\ln\hat{r}/\ln \lambda+\mathrm{i}2m\pi,m\in\mathbb{Z}$. The pre-image of the sequence $B$ is $\{\ln w_n/\ln \lambda+\mathrm{i}2m\pi,m\in\mathbb{Z}\}$. The Riemann surface of $y$ has infinite many sheets. 

We would like to remark that by changing $w\to w^N$ in (\ref{th2.3eq:4}), we obtain
\begin{equation}
\label{eq:7.25}
g(w)=g_1w^N+\sum_{n\geq N}g_nw^{nN},
\end{equation}
which is still a solution of (\ref{th2.3eq:3}) with $b^N=\lambda$. The illustration of the Riemann surface of the continuation of (\ref{eq:7.25}) is the similar to (\ref{th2.3eq:4}), however, instead of having the initial branch at $w_0$, we have $N$ initial branches $w_0\e^{\mathrm{i}2k\pi/N},k=0,\dots,N-1$. When mapping (\ref{eq:7.25}) back into $y(z)$, we still obtain the same solution as when we derived it from (\ref{th2.3eq:4}). By using the uniqueness result from Theorem \ref{auto-lambdasmall}, we can also conclude that all solutions of $g$ of the form $g(w)=\tilde{g}_1w^N+\sum_{n\geq2}\tilde{g}_nw^{n+N}$ must be of the form given in (\ref{eq:7.25}).

\section{The proof of Theorem \ref{th2.2} }\label{sec4}

We will now turn our attention to equation \eqref{eq2.1} with non-constant coefficients by proving Theorem \ref{th2.2}. Let $I$ denote a set of multi-indices $\tau=(j_1,\ldots,j_q)$, where $j_k\in\N\cup\{0\}$ for all $k=1,\ldots,q$. The degree of $\tau$ is defined to be
\begin{equation*}
    d(\tau)=j_1+j_2+\cdots+j_q
\end{equation*}
and the weight of $\tau$ is
	\begin{equation*}
	    w(\tau)=j_1+2j_2+\cdots+qj_q.
	\end{equation*}

We will now transform equation \eqref{eq2.1} into a more convenient form for our purposes. First, by writing \eqref{eq2.1} in the form
	\begin{equation}\label{eq4.1}
	\begin{split}
	y(z+1)&= \left(\sum_{j=1}^{p} a_j(z)y(z)^j\right) \frac{1}{\displaystyle{1+\sum_{k=1}^{q} b_k(z)y(z)^k}}\\
	&= \left(\sum_{j=1}^{p} a_j(z)y(z)^j\right)\sum_{n=0}^{\infty}\left(-\sum_{k=1}^{q} b_k(z)y(z)^k\right)^n\\
	&= \left(\sum_{j=1}^{p} a_j(z)y(z)^j\right)\sum_{n=0}^{\infty}\sum_{\tau\in I \atop d(\tau)=n}(-1)^{d(\tau)} \frac{n!}{j_1! j_2!\cdots j_q!} b_1(z)^{j_1} b_2(z)^{j_2}\cdots b_{q}(z)^{j_{q}}y(z)^{w(\tau)}\\ 
        &= \sum_{j=1}^{p} \sum_{n=0}^{\infty}   \sum_{\tau\in I \atop d(\tau)=n}(-1)^{d(\tau)} \frac{n!}{j_1! j_2!\cdots j_q!}a_j(z)b_1(z)^{j_1} b_2(z)^{j_2}\cdots b_{q}(z)^{j_{q}}y(z)^{w(\tau)+j},
	\end{split}
	\end{equation}
where $a_1(z)\equiv \lambda$, we can see that
	\begin{equation}\label{eq4.2}
	y(z+1)-\lambda y(z)= \sum_{j=2}^{\infty} c_j(z) y(z)^j,
	\end{equation}
where the coefficients $c_j(z)$ are meromorphic functions of $z$, depending on $a_1,\ldots,a_p,b_0,\ldots,b_q$. In particular,
	\begin{equation}\label{eq4.3}
	 \left\{ \begin{array}{ll}
	 c_2(z) &= a_2(z) - a_1(z) b_1(z)    \\
	 c_3(z) &=  a_3(z) - a_2(z) b_1(z) -a_1(z) b_2(z) + a_1(z) b_1(z)^2\\
	 c_4(z) &= a_4(z)-a_3(z)b_1(z)+a_2(z)b_1(z)^2-a_2(z)b_2(z)-
	 a_1(z)b_1(z)^3\\ &\quad +a_1(z)b_1(z)b_2(z)+a_1(z)b_1(z)b_2(z)-a_1(z)b_3(z).\\
	&\,\, \vdots \\
	 \end{array}\right.
	\end{equation}
From \eqref{eq4.3}, we obtain that $c_{2}(z)$ is composed of $2=2^{1}$ terms, $a_2(z)$ and $- a_1(z) b_1(z),$ $c_{3}(z)$ is composed of $4=2^{2}$ terms, and $c_{4}(z)$ is composed of $8=2^{3}$ terms. By an induction argument, it follows that $c_{j}(z)$ is composed of $2=2^{j-1}$ terms. Moreover, the only term in $c_j(z)$ that has degree $j$ is $\pm a_1(z)b_1(z)^{j-1}$.
Hence, we obtain
 \begin{equation}\label{eq2222}
    |c_j(z)| \leq \sum_{k\leq j-1}\max_{i=1,\dots,p}|a_i(z)| \max_{i=1,\dots,p}|b_i(z)|^{k-1}+|a_1||b_1|^{j-1}\leq|\lambda|2^{j-1}\nu^{(j-1)|z|}.
 \end{equation}

Assume that $y(z)$ satisfies \eqref{eq2.4}. Denoting $w(z)=y(z)\lambda^{-z}$, equation \eqref{eq4.2} becomes
	\begin{equation}\label{eq4.4}
	w(z+1)-w(z)=\sum_{j=2}^{\infty} c_j(z) \lambda^{(j-1)z-1}w(z)^j.
	\end{equation}
Therefore,
	\begin{equation}\label{eq4.5}
	 \left\{ \begin{array}{ll}
	 w(z)-w(z-1) &= \displaystyle{\sum_{j=2}^{\infty}} c_j(z-1) \lambda^{(j-1)(z-1)-1}w(z- 1)^j\\
	 w(z-1)-w(z-2) &= \displaystyle{\sum_{j=2}^{\infty}} c_j(z-2) \lambda^{(j-1)(z-2)- 1}w(z-2)^j\\
	   &\,\,\vdots \\
	  w(z-m+1)-w(z-m) &=\displaystyle{ \sum_{j=2}^{\infty}} c_j(z-m) \lambda^{(j-1)(z- m)-1}w(z-
	  m)^j.\\
	 \end{array}\right.
	\end{equation}
By adding up \eqref{eq4.5}, we obtain
	\begin{equation}\label{eq4.6}
	w(z)-w(z-m) =  \sum_{j=2}^{\infty} \sum_{k=1}^m   c_j(z-k)
	\lambda^{(j-1)(z-k)-1} w(z-k)^j.
	\end{equation}
By letting $m\longrightarrow\infty$ in \eqref{eq4.6}, assumption \eqref{eq2.4} yields
	\begin{equation}\label{eq4.7}
	w(z) = \alpha +    \sum_{j=2}^{\infty} \sum_{k=1}^{\infty}  c_j(z-k)
	\lambda^{(j-1)(z-k)-1} w(z-k)^j.
	\end{equation}

We have now transformed equation \eqref{eq2.1} into the desired form. Now, define an operator $T$ by
	\begin{equation}\label{eq4.8}
	T[w](z) = \alpha +   \sum_{j=2}^{\infty} \sum_{k=1}^{\infty}   c_j(z-k)
	\lambda^{(j-1)(z-k)-1} w(z-k)^j.
	\end{equation}
We will define a particular Banach space and show that $T$ is a contraction mapping in this space. 

Note that
	\begin{equation}\label{eq4.9}
	\left|w(z)\right| < L
	\end{equation}
and
	\begin{equation}\label{eq4.10}
	\left|u(z)^j-v(z)^j\right| < jK^{j-1}\left|u(z)-v(z)\right|
	\end{equation}
for suitable constants $L$ and $K$ and for all $(z,w)$, $(z,u)$, $(z,v)$  in the set
	\begin{equation}\label{eq4.11}
 \begin{split}
     \widetilde{D}(\rho,\sigma)=\{(z,w): \Re(z)<-\rho,~|\Im(z)|<\sigma,~|w-\alpha|\leq b\},
 \end{split}
	\end{equation}
where $\rho>0$, $\sigma>0$ and $b$ are fixed constants. Let $X$ be the set of all functions $z\longrightarrow g(z)$, analytic and bounded in
	\begin{equation}\label{eq4.12}
	D(\rho,\sigma)=\left\{(z: \Re(z)< -\rho,,~|\Im(z)|<\sigma, \right\},
	\end{equation}
for which $g(z)\longrightarrow \alpha$ as $\Re(z)\longrightarrow -\infty$ and $\|g-\alpha\|\leq b$,
where
	\begin{equation}\label{eq4.13}
	\left\| g-\alpha \right\| = \sup_{z\in D(\rho,\sigma)}\left| g(z)-\alpha \right|.
	\end{equation}
Clearly $X$ is a metric space under the $\sup$-norm. Moreover, if a sequence $\{g_n\}$ of $X$ is Cauchy, then $\lim_{n\rightarrow\infty} g_n(z)=g(z)$ exists uniformly in $D(\rho,\sigma)$. Furthermore, $g$ is analytic and bounded in $D(\rho,\sigma)$, and
	\begin{equation}\label{eq4.14}
	\lim_{\Re(z)\rightarrow -\infty} g(z) = \lim_{\Re(z)\rightarrow -\infty}
	\lim_{n\rightarrow\infty} g_n(z) = \lim_{n\rightarrow\infty} \lim_{\Re(z)\rightarrow -\infty}
	 g_n(z) = \alpha
	\end{equation}
by uniform convergence. Let finally $\varepsilon >0$. Then there exists $N\in\N$ such that
	\begin{equation}\label{eq4.15}
	\|g-\alpha\| \leq \|g-g_N\| + \|g_N-\alpha \| < b+\varepsilon.
	\end{equation}
Since \eqref{eq4.15} holds for any $\varepsilon >0$,
	\begin{equation}\label{eq4.16}
	\|g-\alpha\|\leq b.
	\end{equation}
Therefore $g\in X$, and so $X$ is a complete metric space.	
	
Next we prove that $T$ is a contraction which maps the space $X$ into itself. For this purpose we need some more information about the set $D(\rho,\sigma)$. Choosing the constant $\rho>0$ in \eqref{eq4.12} sufficiently large, the following four assertions are valid for all $z\in D(\rho,\sigma)$:
	\begin{itemize}
	\item[(i)]  None of the coefficients $c_j$ of \eqref{eq2.1} have poles in $D(\rho,\sigma)$.
	\item[(ii)] Using \eqref{eq2222}, we have for all $j\in \{2,3,\ldots\}$,
\begin{equation}\label{eq3333}
\begin{split}
     |c_j(z)|L^j\leq|\lambda|2^{j-1}L^j\nu^{(j-1)|z|}&\leq|\lambda|^{\epsilon(j-1)(|\Re(z)|+|\Im(z)|)}
     \\
     &=
     |\lambda|^{\epsilon(j-1)(-\Re(z)+|\Im(z)|)},
\end{split}
 \end{equation}
where $0<\epsilon<1$ can be chosen to be fixed and independent of $\rho,$ $\sigma$. This is possible since $\nu<|\lambda|.$
	\item[(iii)] Similarly,
  \begin{equation}\label{eq4444}
     |c_j(z)|jK^{j-1}\leq|\lambda|j2^{j-1}\nu^{(j-1)|z|}\leq|\lambda|^{\epsilon(j-1)(-\Re(z)+|\Im(z)|)}.
 \end{equation}
	\item [(iv)] Let $b$ be the constant in \eqref{eq4.11}. Then 
 \begin{equation}\label{eq4.17}
    \frac{1}{|\lambda|(|\lambda|-1)}\frac{M_S|\lambda|^{(1-\epsilon)\Re(z)}}{1-M_S|\lambda|^{(1-\epsilon)\Re(z)}}< \min\left(b,\tilde{k}\right),
 \end{equation}
	where $\tilde{k}<1$ and $M_S=\sup_{z\in D(\rho,\sigma)}\{|\lambda|^{\epsilon|\Im(z)|-\theta\Im(z)/\ln|\lambda|}\}$.
 	\end{itemize}
We continue by showing that $T[g]$ is analytic in $D(\rho,\sigma)$  whenever $g$ is. Since by \eqref{eq4.9}
\begin{equation}\label{eq5555}
    \begin{split}
        |T[g](z)-\alpha | &= \left|  \sum_{j=2}^{\infty} \sum_{k=1}^{\infty}  c_j(z-k)\lambda^{(j-1)(z-k)-1} w(z-k)^j \right| \\
	& \leq  \sum_{j=2}^{\infty} \sum_{k=1}^{\infty} \left |c_j(z-k)\right|\left|\lambda \right|^{(j-1)(\Re(z)-\theta\Im(z)/\ln|\lambda|-k)-1} \|w^j\|\\
	& \leq \sum_{j=2}^{\infty} \sum_{k=1}^{\infty}\left|c_j(z-k)\right| L^j \left|\lambda \right|^{(j-1)(\Re(z)-\theta\Im(z)/\ln|\lambda|-k)-1}.
    \end{split}
\end{equation}
By using $(ii)$, \eqref{eq5555} becomes	
	\begin{equation}\label{eq4.18}
	\begin{split}
	|T[g](z)-\alpha |
	& \leq \sum_{j=2}^\infty\sum_{k=1}^\infty |\lambda|^{\epsilon(j-1)(-\Re(z)+|\Im(z)|)}|\lambda|^{(j-1)(\Re(z)-\theta\Im(z)/\ln|\lambda|-k)-1} \\
	&= \frac{1}{|\lambda|} \sum_{j=2}^\infty |\lambda|^{(j-1)(1-\epsilon)\Re(z)}M_S^{j-1}\sum_{k=1}\left(|\lambda|^{(j-1)}\right)^{-k} \\
	&= \frac{1}{|\lambda|} \sum_{j=2}^\infty\left( |\lambda| ^{(1-\epsilon)\Re(z)}M_S    \right)^{j-1}\frac{1}{|\lambda|^{j-1}-1} \\
        &\leq \frac{1}{|\lambda|(|\lambda|-1)}\sum_{j=2}^\infty\left( |\lambda| ^{(1-\epsilon)\Re(z)}M_S    \right)^{j-1}\\ 
        &\leq \frac{1}{|\lambda|(|\lambda|-1)}\frac{M_S|\lambda|^{(1-\epsilon)\Re(z)}}{1-M_S|\lambda|^{(1-\epsilon)\Re(z)}}.\\
	\end{split}
	\end{equation}
Therefore the sum $T[g](z)-\alpha$ converges absolutely and uniformly in $D(\rho,\sigma)$, and so $T[g](z)$ is analytic in $D(\rho,\sigma)$. Moreover, by \eqref{eq4.18} and property $(iv)$, $T[g](z)\longrightarrow\alpha$ as $\Re(z)\longrightarrow-\infty$, and
	\begin{equation}\label{eq4.19}
	\| T[g]-\alpha\| \leq b.
	\end{equation}
We conclude that $T$ maps $X$ into itself.

We still need to show that $T$ is a contraction. Let $g,h\in X$. Then by \eqref{eq4.10} and $(iii)$
	\begin{equation}\label{eq4.20}
	\begin{split}
	|T[g](z)-T[h](z) | &= \left| \sum_{j=2}\sum_{k=1} c_j(z-k)\lambda^{(j-1)(z-k)-1} \left(g(z-k)^j-h(z-k)^j\right) \right| \\
	& \leq \sum_{j=2} \sum_{k=1}  \left|c_j(z-k)\right| jK^{j-1} \left|\lambda \right|^{(j-1)(\Re(z)-\theta\Im(z)/\ln|\lambda|-k)-1} \|g-h\|\\
	& \leq \frac{1}{|\lambda|(|\lambda|-1)}\frac{M_S |\lambda|^{(1-\epsilon)\Re(z)}}{1-M_S |\lambda|^{(1-\epsilon)\Re(z)}} \|g-h\|, \\
	\end{split}
	\end{equation}
and so, by property $(iv),$
	\begin{equation}\label{eq4.21}
	\|T[g]-T[h]\| \leq \tilde{k} \|g-h\|,
	\end{equation}
where $\tilde{k}<1$. Hence $T$ is a contraction as asserted.

Since by \eqref{eq4.21} $T$ is a contraction defined on a complete metric space $X$, we conclude by Theorem \ref{th2.2} that $X$ has a unique fixed point under $T$. In other words, there is an analytic function $w\in X$ such that
	\begin{equation}\label{eq4.22}
	w(z) = T[w](z) = \alpha +  \sum_{j=2}^{\infty} \sum_{k=1}^{\infty} c_j(z-k)\lambda^{(j-1)(z-k)-1} w(z-k)^j.
	\end{equation}
Hence, $w$ satisfies equation \eqref{eq4.4}, and so $y(z)=\lambda^z w(z)$ is a solution of \eqref{eq2.1} in $D(\rho,\sigma)$. Finally, using \eqref{eq2.1}, we may continue the solution $y(z)$ into a meromorphic function in the set
\begin{equation*}
    \left\{z\in\C :   |\Im(z)|<\sigma,  \right\}.
\end{equation*}
Since $\sigma>0$ was chosen arbitrarily, the solution is meromorphic in the whole complex plane.

\section{The proof of Theorem \ref{th2.3}}\label{sec5}

The proof of Theorem \ref{th2.3} follows by repeating the reasoning in Section \ref{sec4}, but, instead of the set $D(\rho,\sigma)$, we look at
	\begin{equation*}
	\widehat{D}(\rho,\sigma)=\left\{z: \Re(z) > \rho, \quad |\Im(z)|<\sigma  \right\},
	\end{equation*}
where $R$ is large enough. Since now $|\lambda|<1$, the convergence of
	\begin{equation*}
	\alpha +   \sum_{j=2}^{\infty} \sum_{k=1}^{\infty}   c_j(z-k) 	\lambda^{(j-1)(z-k)-1} w(z-k)^j.
	\end{equation*}
is assured in $\widehat{D}(\rho,\sigma)$. In the following proof, it is only necessary to make the following changes to equations \eqref{eq4.17}, \eqref{eq4.18} and \eqref{eq4.20} in Section \ref{sec4}:
\begin{enumerate}
    \item [(1)] Replace \eqref{eq4.17} with 
    \begin{equation*}
         \frac{1}{|\lambda|(1-|\lambda|)}\frac{|\lambda|^{(1-\epsilon)\Re(z)}M_S}{1-|\lambda|^{(1-\epsilon)\Re(z)}M_S}< \min\left(b,\tilde{k}\right),
    \end{equation*}
    \item [(2)] Replace \eqref{eq4.18} with 
    \begin{equation*}
        |T[g](z)-\alpha |\leq\frac{1}{|\lambda|(1-|\lambda|)} \frac{|\lambda|^{(1-\epsilon)\Re(z)}M_S}{1-|\lambda|^{(1-\epsilon)\Re(z)}M_S}.
    \end{equation*}
    \item [(2)] Replace \eqref{eq4.20} with 
    \begin{equation*}
        |T[g](z)-T[h](z) |\leq \frac{1}{|\lambda|(1-|\lambda|)}\frac{|\lambda|^{(1-\epsilon)\Re(z)}M_S}{1-|\lambda|^{(1-\epsilon)\Re(z)}M_S}\|g-h\|.
    \end{equation*}
\end{enumerate}
Then the rest of the proof of the local existence and uniqueness of the solution in $\widehat{D}(\rho,\sigma)$ is identical to what was done in Section \ref{sec4}, and we will not repeat it here. 

\section{Discussion}
In this paper we have studied the existence of solutions of first-order difference equations in the complex domain.  We proved the existence of meromorphic solutions for a class of non-autonomous equations with growth conditions on the coefficients as $\Re(z)\to-\infty$. We also provided details of a direct proof of a known result on the existence of meromorphic solutions in the constant coefficient case in which the multiplier at a fixed point of $R$ is one.  Additionally, we conducted a detailed analysis of the infinite-sheeted Riemann surface on which the analytic continuation of a solution to a simple constant coefficient equation is defined.

\bibliographystyle{amsplain}
\bibliography{paper2}

\def\cprime{$'$}
\providecommand{\bysame}{\leavevmode\hbox to3em{\hrulefill}\thinspace}
\providecommand{\MR}{\relax\ifhmode\unskip\space\fi MR }
\providecommand{\MRhref}[2]{%
  \href{http://www.ams.org/mathscinet-getitem?mr=#1}{#2}
}
\providecommand{\href}[2]{#2}
\begin{thebibliography}{10}

\bibitem{p2Azarina}
Y.~V. Azarina, \emph{Meromorphic solutions of the equation
  {{\(\omega(z+1)=R[\omega(z)]\)}}}, Teor. Funkts. Funkts. Anal. Prilozh.
  \textbf{48} (1987), 26--32 (Russian).

\bibitem{p2BakerLiverpool}
I.~N. Baker and L.~S.~O. Liverpool, \emph{The entire solutions of a polynomial
  difference equation}, Aequationes Math. \textbf{27} (1984), 97--113.

\bibitem{fatou1919}
P.~Fatou, \emph{Sur les \'{e}quations fonctionnelles}, Bull. Soc. Math. France
  \textbf{47} (1919), 161--271.

\bibitem{fatou1920a}
\bysame, \emph{Sur les \'{e}quations fonctionnelles}, Bull. Soc. Math. France
  \textbf{48} (1920), 33--94.

\bibitem{fatou1920b}
\bysame, \emph{Sur les \'{e}quations fonctionnelles}, Bull. Soc. Math. France
  \textbf{48} (1920), 208--314.

\bibitem{p2HKphD2006}
R.~G. Halburd and R.~J. Korhonen, \emph{Existence of finite-order meromorphic
  solutions as a detector of integrability in difference equations}, Phys. D
  \textbf{218} (2006), no.~2, 191--203.

\bibitem{p2julia}
G.~Julia, \emph{Memoire sur l'iteration des fonctions rationnelles}, Journal de
  Mathematiques Pures et Appliquees \textbf{1} (1918), 47--246 (French).

\bibitem{p2kimura1971}
T.~Kimura, \emph{On the iteration of analytic functions}, Funkcial. Ekvac.
  \textbf{4} (1971), 197--238.

\bibitem{p2kimura1973}
\bysame, \emph{On meromorphic solutions of the difference equation
  {$y(x+1)=y(x)+1+(\lambda /y(x) )$}}, Symposium on {O}rdinary {D}ifferential
  {E}quations ({U}niv. {M}innesota, {M}inneapolis, {M}inn., 1972; dedicated to
  {H}ugh {L}. {T}urrittin), Lecture Notes in Math., Vol. 312, Springer, Berlin,
  1973, pp.~74--86.

\bibitem{koenigs:84}
G.~Koenigs, \emph{Recherches sur les int\'{e}grales de certaines \'{e}quations
  fonctionnelles}, Ann. Sci. \'{E}cole Norm. Sup. (3) \textbf{1} (1884), 3--41.

\bibitem{mahler1981special}
K.~Mahler, \emph{On a special nonlinear functional equation}, Proceedings of
  the Royal Society of London. A. Mathematical and Physical Sciences
  \textbf{378} (1981), no.~1773, 155--178.

\bibitem{poincare:90}
G.~Poincar\'e, \emph{Sur une classe nouvelle de transcendantes uniformes}, J.
  Math. Pures Appl. \textbf{6} (1890), 313--365.

\bibitem{p2Simo}
S.~Shimomura, \emph{Entire solutions of a polynomial difference equation}, J.
  Fac. Sci. Univ. Tokyo Sect. IA Math. \textbf{28} (1981), no.~2, 253--266.

\bibitem{wedderburn1922mahler}
J.~H.~M. Wedderburn, \emph{The functional equation $g(x^2) = 2ax + [g(x)]^2$},
  Ann. of Math. \textbf{24} (1922), 121--140.

\bibitem{p2Yana}
N.~Yanagihara, \emph{Meromorphic solutions of some difference equations},
  Funkcial. Ekvac. \textbf{23} (1980), no.~3, 309--326.

\end{thebibliography}

\end{document}